\documentclass[12pt]{article}
\usepackage[margin=1in]{geometry}
\usepackage{graphicx}    
\usepackage{amsmath} 
\usepackage{xcolor}
\usepackage{amssymb}
\usepackage{amsthm} 
\usepackage{amsfonts}
\usepackage{booktabs}
\usepackage{hyperref}    
\usepackage{geometry}    
\usepackage{caption}     
\usepackage{cite}
\usepackage{algorithm}
\usepackage{algorithmicx}
\usepackage{algpseudocode}
\usepackage{pgfplots}
\usepackage{bm}
\usepackage{float}
\usetikzlibrary{patterns,positioning}
\usepackage{tikz}
\usetikzlibrary{decorations.markings}
\usetikzlibrary{matrix}
\usepackage{pgfplots}
\usetikzlibrary{shapes, arrows.meta, patterns}
\usepackage{mathrsfs}
\usepackage{siunitx}

\usepackage{subcaption}

\newcommand{\bfa}[1]{\boldsymbol{#1}}

\newcommand{\Sym}{\text{Sym}}   			%

\DeclareMathAlphabet{\mathpzc}{OT1}{pzc}{m}{it}

\usepackage{authblk}

\usepackage{etoolbox}

\newtheorem{remark}{Remark}

\providecommand{\keywords}[1]
{
  \small	
  \textbf{\textit{Keywords---}} #1
}

\title{An HHT-$\alpha$-based finite element framework for wave propagation in constitutively nonlinear elastic materials}

\author{S. M. Mallikarjunaiah}

\affil{Department of Mathematics \& Statistics, Texas A\&M University-Corpus Christi, TX- 78412, USA}
\affil[ ]{\textit{E-mail addresses:} \texttt{M.Muddamallappa@tamucc.edu}}
\date{}

\begin{document}

\maketitle  

\begin{abstract}
This paper presents a computational framework for modeling wave propagation in geometrically linear elastic materials characterized by algebraically nonlinear constitutive relations. We derive a specific form of the nonlinear wave equation in which the nonlinearity explicitly appears in the time-derivative terms that govern the evolution of the mechanical fields. The numerical solution is established using a fully discrete formulation that combines the standard finite element method for spatial discretization with the implicit \textit{Hilber-Hughes-Taylor} (HHT)-$\alpha$ scheme for time integration. To address the nonlinear nature of the discrete system, we employ Newton’s method to iteratively solve the linearized equations at each time step. The accuracy and robustness of the proposed framework are rigorously verified through convergence analyses, which demonstrate optimal convergence rates in both space and time. Furthermore, a detailed parametric study is conducted to elucidate the influence of the model's constitutive parameters. The results reveal that the magnitude parameter of the stress-dependent variation in wave speed leads to wavefront steepening and the formation of shock discontinuities. Conversely, the exponent parameter acts as a nonlinearity filter; high values suppress nonlinear effects in small-strain regimes, whereas low values allow significant dispersive behavior. This work provides a validated tool for analyzing shock formation in advanced nonlinear materials.
\end{abstract}

\vspace{.1in}

\noindent \keywords{finite element method; Nonlinear material; HHT-$\alpha$ time discretization; Nonlinear hyperbolic PDE; Wave speed}
 
\section{Introduction}

The accurate computational modeling of wave propagation in solid media is ubiquitous across engineering disciplines, serving as a fundamental tool for both safety assessment and technological innovation. In the aerospace and automotive sectors, it underpins modern Non-Destructive Testing techniques, which are essential for detecting microscopic defects in composite materials before catastrophic failures occur \cite{Rose2023NDT, Giurgiutiu2024Composites}. Similarly, in civil engineering and geophysics, high-fidelity wave simulations are imperative for quantifying seismic risks and designing infrastructure capable of withstanding dynamic shock loads \cite{Bozorgnia2023Seismic, Chopra2024Dynamics}. Beyond structural mechanics, these models are increasingly vital in biomedical engineering, particularly for refining ultrasound elastography methods used in early cancer detection \cite{Parker2023Elastography}. As industries transition toward advanced materials with complex nonlinear behaviors, the ability to robustly predict dispersive phenomena and shock formation becomes critical for ensuring operational reliability in these high-stakes real-world applications \cite{LeVeque2002}.

The conventional modeling of wave propagation in elastic solids predominantly relies on linear elastodynamic theory, which assumes a constant linear relationship between stress and strain. While computationally convenient, this linearization is often unphysical for describing the dynamic response of many real-world materials, as it inherently enforces a globally constant wave speed $c_0 = \sqrt{E/\rho}$ regardless of the deformation state \cite{Achenbach1973}. In contrast, a wide range of media—including polymers, biological soft tissues, and granular soils—exhibit strong constitutive nonlinearity, where the instantaneous stiffness, and consequently the wave speed, evolves as a function of the local stress amplitude \cite{Destrade2011}. By neglecting this stress-velocity coupling, linear theory fails to capture critical nonlinear phenomena such as wavefront steepening and the eventual formation of shock discontinuities, which are the dominant mechanisms in the high-energy dynamics of physically nonlinear materials \cite{LeVeque2002}.

The mathematical modeling of complex material behavior has been significantly broadened by the introduction of implicit constitutive theories, a framework extensively developed by Rajagopal and co-workers. Traditionally, the Cauchy theory of elasticity presumes an explicit functional relationship where the stress tensor $\mathbf{T}$ is determined solely by the deformation gradient or the linearized strain $\boldsymbol{\varepsilon}$, i.e., $\mathbf{T} = \mathcal{G}(\boldsymbol{\varepsilon})$. However, Rajagopal \cite{Rajagopal2003, Rajagopal2006} argued that this explicit definition is overly restrictive and fails to capture the response of materials that do not exhibit a one-to-one correspondence between stress and strain. By generalizing the constitutive structure to an implicit relation of the form $\boldsymbol{f}(\mathbf{T}, \boldsymbol{\varepsilon}) = \mathbf{0}$, Rajagopal established a unified theoretical basis that encompasses classical elasticity, viscoelasticity, and plasticity as special cases. This implicit framework is thermodynamically consistent \cite{RajagopalSrinivasa2007} and is particularly powerful because it allows the kinematical variables to be expressed as functions of the stress, rather than the reverse. This inversion is crucial for describing biological tissues, polymers, and other ``physically nonlinear'' materials where the compliance changes drastically with the applied load, a behavior that standard explicit models struggle to represent accurately.

A critical subclass of these implicit materials, and the focus of the present work, is those governed by algebraically nonlinear relations within the small-strain regime. While classical linearized elasticity assumes a constant fourth-order stiffness tensor, Rajagopal's implicit theory permits the material moduli to depend on the stress state itself, even when geometric deformations remain infinitesimal \cite{Rajagopal2011}. This leads to the concept of ``limiting strain'' or ``limiting stress'' models, where the material response remains bounded despite increasing loads. Such models are instrumental in resolving physical paradoxes found in linear theories, such as the prediction of infinite strains at stress concentration points. As demonstrated by Bustamante and Rajagopal \cite{Bustamante2011}, these algebraically nonlinear relations describe a material wherein the wave speed is not a constant $c_0$ but a function of the local stress field. This stress-dependency is precisely the mechanism that drives the dispersive phenomena and shock formation observed in the numerical simulations presented in this study, confirming that physical nonlinearity can exist independently of geometric nonlinearity.

The implications of implicit constitutive relations extend significantly into the field of fracture mechanics, where they offer a regularization of the singularities inherent to Linear Elastic Fracture Mechanics (LEFM). In standard LEFM, the stress and strain fields are predicted to become singular at the crack tip ($r^{-1/2}$ singularity), a physical impossibility. However, by utilizing Rajagopal's nonlinear implicit models, specifically those with strain-limiting characteristics, it is possible to bound the strain energy density and deformation at the crack tip, providing a more realistic description of the near-tip fields. Recent computational efforts have successfully incorporated these models into robust numerical frameworks. Notably, the author has previously utilized this approach to develop finite element framework for static cracks, demonstrating that implicit relations can be effectively discretized to capture complex fracture topologies \cite{MP2025hp,mallikarjunaiah2025hp,gou2023computational,gou2023finite,mallikarjunaiah2025crack,mallikarjunaiah2015}. Furthermore, this constitutive class has been extended to phase-field modeling of fracture, where the nonlinearity in strain-energy density was utilized in udnerstanding the propagation of cracks \cite{manohar2025adaptive,fernando2025xi,fernando2025,manohar2025convergence}.

The investigation of wave propagation in materials governed by implicit constitutive relations has gained significant traction following the foundational work of Rajagopal \cite{Rajagopal2003, Rajagopal2006}, who generalized the standard Cauchy elasticity to include stress-dependent moduli within a thermodynamically consistent framework. Unlike classical linear theories where wave speed is invariant, these implicit models naturally predict a nonlinear dependence of the wave velocity on the local stress state, a phenomenon critically analyzed in many works \cite{bustamante2015direct,Bustamante2011,meneses2018note,ibarra2022analysis,bustamante2024iteration,kannan2014unsteady}. Their studies demonstrated that even in the absence of geometric nonlinearity, the constitutive nonlinearity inherent to these ``physically nonlinear'' bodies can drive the formation of shock discontinuities from smooth initial conditions. Further extending this analysis, Magan et al. \cite{Magan2019} investigated power-law implicit materials in cylindrical geometries, confirming that the steepening of wavefronts and the subsequent shock development are governed by the specific exponents of the constitutive relation. Recent contributions by Bustamante et al. \cite{Bustamante2025} have expanded these concepts to stretch-limited elastic strings, revealing unique longitudinal shock behaviors that deviate significantly from traditional hyperelastic predictions. Collectively, this body of literature establishes that implicit constitutive theories provide a rich and necessary platform for capturing the complex, dispersive dynamics of modern advanced materials.

The primary contribution of this study is the development of a robust computational framework for simulating transient wave propagation in geometrically linear, physically nonlinear elastic solids. A key theoretical advancement presented herein is the derivation of a second-order nonlinear wave equation governing the stress field, formulated such that the nonlinearity resides explicitly within the time-dependent terms. While pure stress formulations are well-established for linear elastodynamics \cite{ignaczak1963rayleigh, norris2021stress}, extensions to nonlinear constitutive behavior typically rely on mixed velocity-stress formulations \cite{abramov2011nonlinear}. To the authors' knowledge, the direct solution of the pure nonlinear stress wave equation derived in this work has not been widely addressed in the computational mechanics literature. To solve this system, we propose a fully discrete numerical scheme that combines a standard continuous Galerkin Finite Element Method (FEM) for spatial discretization with the implicit HHT-$\alpha$ method for time integration. The HHT-$\alpha$ algorithm is specifically chosen for its second-order accuracy and its ability to provide controllable numerical dissipation, effectively damping high-frequency spurious oscillations. The resulting nonlinear discrete equations are solved iteratively at each time step using a Newton-Raphson procedure. Finally, we conduct a systematic parametric study to elucidate the roles of the constitutive parameters, characterizing the transition from non-dispersive linear wave propagation to highly dispersive, shock-dominated dynamics.

The remainder of this paper is organized as follows. Section \ref{sec:formulation} establishes the theoretical framework, detailing the class of implicit constitutive relations considered and deriving the governing nonlinear wave equation for the stress field. Section \ref{sec:numerical} presents the fully discrete numerical formulation, combining the continuous Galerkin finite element method for spatial discretization with the HHT-$\alpha$ time-integration scheme and the Newton-Raphson linearization procedure. Section \ref{sec:results} discusses the numerical results, providing a comprehensive verification against linear theory and a detailed parametric study on the influence of the constitutive parameters $b$ and $a$ on the wave propagation in physically nonlinear solids. Finally, Section \ref{sec:conclusion} summarizes the key findings and offers concluding remarks.

\section{Mathematical formulation of Rajagopal's implicit theory of elasticity}\label{sec:formulation}

\subsection{Kinematics}
Let $\mathcal{B}$ be a continuum body identified with its reference configuration $\Omega_R \subset \mathbb{R}^d$ ($d=1,2,3$). The motion is defined by a smooth mapping $\boldsymbol{\chi}: \Omega_R \times [0, T] \to \mathbb{R}^d$, such that the current position is $\mathbf{x} = \boldsymbol{\chi}(\mathbf{X}, t)$. The local deformation is characterized by the deformation gradient $\mathbf{F} = \nabla_{\mathbf{X}} \boldsymbol{\chi}$, which admits the polar decomposition $\mathbf{F} = \mathbf{R}\mathbf{U} = \mathbf{V}\mathbf{R}$, where $\mathbf{R} \in \text{SO}(d)$ is the rotation, and $\mathbf{U}, \mathbf{V} \in \text{Sym}^+(d)$ are the right and left stretch tensors, respectively. The associated Cauchy-Green tensors are defined as $\mathbf{C} = \mathbf{F}^T \mathbf{F}$ and $\mathbf{B} = \mathbf{F}\mathbf{F}^T$.

We restrict our analysis to the regime of \textit{small displacement gradients}. Assuming $\max_{\mathbf{x},t} |\nabla \mathbf{u}| = \mathcal{O}(\delta)$ with $\delta \ll 1$, the Green-Lagrange strain reduces to the standard linearized strain tensor $\boldsymbol{\epsilon}$:
\begin{equation}
    \boldsymbol{\epsilon} = \frac{1}{2} [\nabla \mathbf{u} + (\nabla \mathbf{u})^T].
\end{equation}

\subsection{Implicit constitutive theory}
Classical Cauchy elasticity posits an explicit relation $\mathbf{T} = \widehat{\mathbf{G}}(\mathbf{B})$, where $\mathbf{T}$ denotes the Cauchy stress. While sufficient for many materials, this explicit framework is restrictive; it cannot capture behaviors such as limiting strains in the small-strain regime \cite{rajagopal2003implicit,rajagopal2007elasticity}. To overcome this, Rajagopal proposed a generalization based on implicit constitutive relations of the form $\mathcal{G}(\mathbf{T}, \mathbf{B}) = 0$.

For isotropic materials, representation theorems for tensor functions \cite{truesdell1992first,gou2015modeling} allow the implicit relation to be expanded generally as:
\begin{equation}
    \chi_0 \mathbf{I} + \chi_1 \mathbf{T} + \chi_2 \mathbf{B} + \dots + \chi_8 (\mathbf{T}^2 \mathbf{B}^2 + \mathbf{B}^2 \mathbf{T}^2) = \mathbf{0},
    \label{eq:general_implicit}
\end{equation}
where the coefficients $\chi_i$ depend on the invariants of $\mathbf{T}$ and $\mathbf{B}$.

A critical divergence from classical theory occurs upon linearization. While linearizing the classical explicit model yields Hooke's Law ($\mathbf{T} = \mathbb{C}:\boldsymbol{\epsilon}$), linearizing the implicit relation \eqref{eq:general_implicit} under the assumption $\mathbf{B} \approx \mathbf{I} + 2\boldsymbol{\epsilon}$ yields a \textbf{nonlinear relationship between linearized strain and stress}:
\begin{equation}
    \boldsymbol{\epsilon} = \hat{\chi}_0 \mathbf{I} + \hat{\chi}_1 \mathbf{T} + \hat{\chi}_2 \mathbf{T}^2 + \mathcal{O}(\delta^2).
\end{equation}
This framework allows for the modeling of materials where strains remain infinitesimal even as stresses become arbitrarily large, a phenomenon known as \textit{strain-limiting behavior}.

\subsection{A specific strain-limiting class}
We adopt a specific subclass of these materials characterized by a compliance-based response function $\bfa{F}: \Sym(d) \to \Sym(d)$ \cite{itou2018states,bulivcek2015analysis}:
\begin{equation}
    \boldsymbol{\epsilon} = \bfa{F}(\bfa{T}) := \frac{\mathbb{K}[\bfa{T}]}{\left(1 + b^{a} \| \mathbb{K}^{1/2}[\bfa{T}] \|^{a} \right)^{1/a}}.
    \label{eq:specific_F}
\end{equation}
Here, $\mathbb{K} = \mathbb{E}^{-1}$ is the fourth-order compliance tensor, and $a, b \geq 0$ are material parameters. The parameter $\beta$ dictates the limiting strain threshold; specifically, the model enforces the uniform bound:
\begin{equation}
    \sup_{\bfa{T} \in \Sym} \| \bfa{F}(\bfa{T}) \| \leq \frac{1}{\beta}.
\end{equation}
This formulation ensures physical plausibility by precluding the singular infinite strains predicted by linear elasticity at crack tips \cite{yoon2022CNSNS}.

\textbf{Mathematical properties.} The response function $\bfa{F}$ defined in \eqref{eq:specific_F} possesses key properties ensuring the well-posedness of the boundary value problem:
\begin{itemize}
    \item[(i)] \textbf{Boundedness:} $\| \bfa{F}(\bfa{T}) \| \leq 1/\beta$ for all $\bfa{T}$, ensuring finite strains globally.
    \item[(ii)] \textbf{Strict Monotonicity:} $(\bfa{F}(\bfa{T}_1) - \bfa{F}(\bfa{T}_2)) \colon (\bfa{T}_1 - \bfa{T}_2) > 0$ for distinct $\bfa{T}_1, \bfa{T}_2$, guaranteeing uniqueness.
    \item[(iii)] \textbf{Lipschitz Continuity:} $\| \bfa{F}(\bfa{T}_1) - \bfa{F}(\bfa{T}_2) \| \leq L \| \bfa{T}_1 - \bfa{T}_2 \|$, ensuring stability.
    \item[(iv)] \textbf{Coercivity:} The operator satisfies coercivity conditions required for existence in the variational setting.
\end{itemize}

\subsection{Application to biomechanics: parameter identification}

To validate the descriptive capability of the proposed model, we calibrate the constitutive parameters against experimental data for biological tissues \cite{sengul2021viscoelasticity}. The 1D reduction of \eqref{eq:specific_F} yields the scalar relation:
\begin{equation}
    \epsilon = \frac{\sigma}{\left(1 + (b |\sigma|)^a\right)^{1/a}}.
    \label{eq:strain_limiting_1D}
\end{equation}
Parameters $a$ and $b$ were identified via nonlinear least-squares optimization to minimize the sum of squared errors (SSE) between the model and experimental data for Thoracic and Carotid artery tissues. The results, summarized in Table \ref{tab:fit_results}, demonstrate a high goodness-of-fit ($R^2 > 0.94$), confirming the model's ability to capture the nonlinear saturation inherent in soft tissue mechanics.

\begin{table}[h!]
    \centering
    \caption{Optimized material parameters and goodness of fit statistics.}
    \label{tab:fit_results}
    \vspace{0.2cm}
    \begin{tabular}{lccc}
        \toprule
        \textbf{Dataset} & $b$ & $a$ & $R^2$ \\
        \midrule
        Thoracic A. (Longitudinal) & $3.8106 \times 10^{-1}$ & $0.1765$ & $0.9780$ \\
        Carotid Artery             & $1.1389 \times 10^{-6}$ & $0.0505$ & $0.9635$ \\
        Thoracic A. (Transverse)   & $2.7353 \times 10^{-8}$ & $0.0433$ & $0.9476$ \\
        \bottomrule
    \end{tabular}
\end{table}

\begin{figure}[H]
    \centering
    \includegraphics[width=0.5\linewidth]{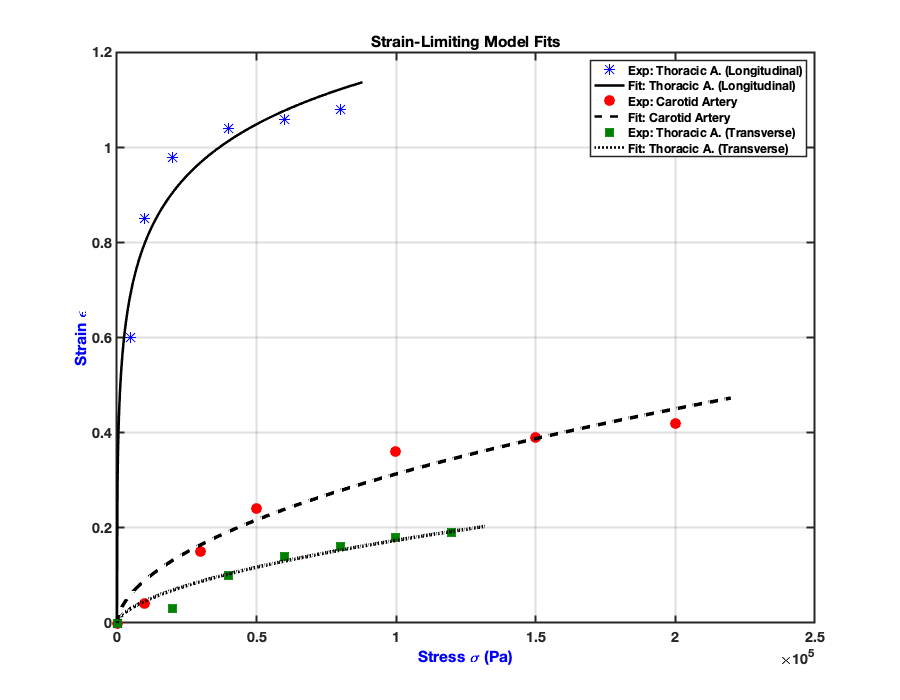}
    \caption{Comparison of experimental stress-strain data (markers) with the implicit strain-limiting model fits (lines). The model accurately captures the saturation behavior of the tissue.}
    \label{fig:strain_limiting_fits}
\end{figure}

\subsection{Derivation of the 1D ``stress'' equation of motion}
To derive the equation of motion for a one-dimensional elastic body defined by the implicit constitutive relation 
\begin{equation}\label{1Dmodel}
\epsilon = \mathfrak{f}(\sigma):=  \frac{\sigma}{\left(1 + (b |\sigma|)^a\right)^{1/a}}.
\end{equation}
Here $\epsilon$ and $\sigma$ are the linearized strain and Cauchy stress in 1D. 
We begin with the standard balance of linear momentum in one dimension. Assuming a constant mass density $\rho$ and a body force $b$, the equation relating the stress $\sigma$ and displacement $u$ is given by the balance of linear momentum:
\begin{equation}
    \frac{\partial \sigma}{\partial x} + \rho b = \rho \frac{\partial^2 u}{\partial t^2}
\end{equation}
To express this solely in terms of stress, we differentiate the entire equation with respect to the spatial coordinate $x$. This allows us to relate the acceleration term directly to the strain $\epsilon$:
\begin{equation}
    \frac{\partial^2 \sigma}{\partial x^2} + \rho \frac{\partial b}{\partial x} = \rho \frac{\partial^2}{\partial t^2} \left( \frac{\partial u}{\partial x} \right)
\end{equation}
In the one-dimensional case, the linearized strain is defined as the spatial gradient of the displacement, 
\begin{equation}
\epsilon = \frac{\partial u}{\partial x}. 
\end{equation}
Substituting this into the right-hand side of the previous equation, we get:
\begin{equation}
    \frac{\partial^2 \sigma}{\partial x^2} + \rho \frac{\partial b}{\partial x} = \rho \frac{\partial^2 \epsilon}{\partial t^2}
\end{equation}
For this specific class of elastic bodies, the strain is given as a nonlinear function of the stress as in \eqref{1Dmodel}. By substituting this constitutive relation into our equation, we arrive at a single partial differential equation governing the stress evolution:
\begin{equation}
    \frac{\partial^2 \sigma}{\partial x^2} + \rho \frac{\partial b}{\partial x} = \rho \frac{\partial^2}{\partial t^2} [\mathfrak{f}(\sigma)]
\end{equation}
If we assume the body forces are negligible ($b=0$), this simplifies to the wave equation that is nonlinear in time but linear in space part:
\begin{equation}\label{stress_wave_eq}
    \frac{\partial^2 \sigma}{\partial x^2} = \rho \frac{\partial^2}{\partial t^2} [\mathfrak{f}(\sigma)]
\end{equation}
\begin{remark}[Stress-based formulation vs. standard elastodynamics]
The partial differential equation derived in Eq.~\eqref{stress_wave_eq} differs fundamentally from the classical displacement-based formulation of nonlinear elastodynamics. In standard hyperelasticity, the equation of motion is typically expressed as $\rho \ddot{u} = \nabla \cdot \sigma(\nabla u)$, resulting in a quasilinear PDE where the nonlinearity resides in the spatial operator (the divergence of the stress). Conversely, the present formulation is derived for the stress variable $\sigma$ directly. Here, the spatial operator $\partial_{xx} \sigma$ remains linear, while the nonlinearity is encapsulated entirely within the temporal inertial term $\partial_{tt} [\mathfrak{f}(\sigma)]$. This structure implies that the effective wave speed is governed by the tangent compliance of the material, which evolves dynamically with the stress state, rather than the tangent stiffness.
\end{remark} 

The governing partial differential equation derived in \eqref{stress_wave_eq} exhibits a unique structure where the nonlinearity is entirely confined to the inertial term. To understand the physical implications, we apply the chain rule to the time derivative of the constitutive function $\mathfrak{f}(\sigma)$. Letting $(\cdot)'$ denote differentiation with respect to stress and $(\dot{\cdot})$ denote differentiation with respect to time, the second time derivative expands as:
\begin{equation}
    \frac{\partial^2}{\partial t^2} [\mathfrak{f}(\sigma)] = \mathfrak{f}'(\sigma) \ddot{\sigma} + \mathfrak{f}''(\sigma) (\dot{\sigma})^2.
\end{equation}
Substituting this expansion back into the stress wave equation yields the quasi-linear form:
\begin{equation} \label{expanded_wave_eq}
    \frac{\partial^2 \sigma}{\partial x^2} = \rho \mathfrak{f}'(\sigma) \frac{\partial^2 \sigma}{\partial t^2} + \rho \mathfrak{f}''(\sigma) \left(\frac{\partial \sigma}{\partial t}\right)^2.
\end{equation}
This expanded form reveals two critical physical characteristics of the proposed formulation:

\begin{itemize}
    \item \textbf{State-dependent wave speed:} The coefficient $\rho \mathfrak{f}'(\sigma)$ acts as an effective, dynamic density. The propagation speed of stress waves, denoted by $c(\sigma)$, is determined by the tangent compliance of the material:
    \begin{equation}
        c(\sigma) = \sqrt{\frac{1}{\rho \mathfrak{f}'(\sigma)}}.
    \end{equation}
    Consequently, the wave speed is not constant but evolves with the local stress state. Regions of the material undergoing stiffening (decreasing compliance $\mathfrak{f}'$) will propagate waves faster, while softening regions will retard wave propagation.

    \item \textbf{Quadratic velocity nonlinearity:} The term $\rho \mathfrak{f}''(\sigma) (\dot{\sigma})^2$ acts as a nonlinear source term driven by the square of the stress rate. This implies that even for materials with moderate static nonlinearity, high-frequency loading (large $\dot{\sigma}$) will generate significant nonlinear feedback. This term is analogous to convective forces in fluids but operates in the temporal domain, acting as an inertial driver that pumps energy into higher harmonics of the solution.
\end{itemize}

\begin{remark}[Well-posedness and stability conditions]
The hyperbolic character of the governing equation depends critically on the monotonicity of the constitutive relation. Rearranging the expanded form in Eq.~\eqref{expanded_wave_eq}, the equation can be viewed as a wave equation with a variable coefficient $\rho \mathfrak{f}'(\sigma)$ scaling the acceleration term $\ddot{\sigma}$. For the problem to remain well-posed (i.e., strictly hyperbolic with real-valued wave speeds), this coefficient must remain strictly positive:
\begin{equation} \label{stability_cond}
    \rho \mathfrak{f}'(\sigma) > 0 \quad \implies \quad \mathfrak{f}'(\sigma) > 0.
\end{equation}
Physically, this condition requires that the tangent compliance of the material remains positive—meaning the material must not exhibit a snap-back instability or negative stiffness in the stress-strain response. If $\mathfrak{f}'(\sigma) \to 0$, the wave speed diverges to infinity, and if $\mathfrak{f}'(\sigma) < 0$, the equation changes type from hyperbolic to elliptic, leading to ill-posedness where the initial value problem is unstable (Hadamard instability). Therefore, the material parameters $a$ and $b$ in the specific constitutive model \eqref{1Dmodel} must be chosen such that the monotonicity of $\mathfrak{f}(\sigma)$ is preserved over the entire range of admissible stress values.
\end{remark}

\subsection{Problem formulation}

\subsubsection{Strong form}
We consider a one-dimensional nonlinear wave equation where the Cauchy stress $\sigma(x,t)$ serves as the primary unknown variable. The domain is defined as $\Omega = (0, L)$ with a time interval $(0, T]$. The governing equation in the strong form is derived from the balance of linear momentum:
\begin{equation} \label{strong_form_general}
    \rho \frac{\partial^2 \epsilon}{\partial t^2} - \frac{\partial^2 \sigma}{\partial x^2} = f(x,t) \quad \text{in } \Omega \times (0, T],
\end{equation}
where $\rho$ is the constant mass density, $f(x,t)$ is the body force, and $\epsilon$ is the linearized strain. The material behavior is governed by the implicit constitutive relation parameterized by material constants $a$ and $b$:
\begin{equation} \label{constitutive_model_ab}
    \epsilon(\sigma) = \frac{\sigma}{\left(1 + (b |\sigma|)^a \right)^{1/a}}.
\end{equation}
To express the momentum balance purely in terms of stress, we apply the chain rule to the acceleration term, expanding the second time derivative of strain as:
\begin{equation}
    \ddot{\epsilon} = \epsilon'(\sigma) \ddot{\sigma} + \epsilon''(\sigma) (\dot{\sigma})^2.
\end{equation}
Substituting this expansion into \eqref{strong_form_general} yields the final strong form of the stress wave equation:
\begin{equation} \label{strong_form_expanded}
    \rho \left[ \epsilon'(\sigma) \ddot{\sigma} + \epsilon''(\sigma) (\dot{\sigma})^2 \right] - \frac{\partial^2 \sigma}{\partial x^2} = f(x,t).
\end{equation}

\subsubsection{Initial and boundary conditions}
To ensure the problem is well-posed, the partial differential equation is closed with the following conditions:

\paragraph{Initial conditions:}
The system is initialized with a known stress distribution and stress rate:
\begin{equation}
    \sigma(x,0) = \sigma_0(x), \quad \dot{\sigma}(x,0) = v_0(x).
\end{equation}

\paragraph{Boundary conditions:}
We apply Dirichlet-type boundary conditions on the stress variable, which correspond to specific physical loading constraints:
\begin{align}
    \sigma(0,t) &= 0, && \text{(Left: Free boundary / Stress-free)} \\
    \sigma(L,t) &= A \sin(\omega t), && \text{(Right: Applied oscillatory traction)}.
\end{align}
\begin{remark}
    Note that a Dirichlet condition on stress ($\sigma = 0$) corresponds physically to a free end. Conversely, a fixed displacement boundary condition (clamped, $u=0$) would manifest in this formulation as a Neumann condition on the stress gradient ($\frac{\partial \sigma}{\partial x} = 0$).
\end{remark}

The initial-boundary value problem defined by Eq.~\eqref{strong_form_expanded} presents specific challenges for numerical solution. While the spatial operator ($\partial_{xx} \sigma$) is linear, the temporal term contains high-order nonlinearities ($\epsilon' \ddot{\sigma}$ and $\epsilon'' \dot{\sigma}^2$) that can act as sources of numerical instability. Standard explicit time-integration schemes may struggle with stability due to the variable "mass" coefficient $\rho \epsilon'(\sigma)$.

To address this, we adopt a semi-discrete approach. We first effectively decouple the spatial and temporal discretizations:
\begin{enumerate}
    \item \textbf{Spatial discretization:} We employ the Galerkin Finite Element Method (FEM) to handle the spatial domain. This allows us to convert the linear Laplacian operator into a constant stiffness matrix, simplifying the system structure.
    \item \textbf{Temporal discretization:} We utilize the Hilbert-Hughes-Taylor (HHT-$\alpha$) implicit time-integration scheme \cite{hilbert1977improved,hughes2012finite}. The HHT-$\alpha$ method is chosen specifically for its ability to introduce controllable algorithmic damping in the high-frequency modes—which are often excited artificially by the nonlinear velocity term $(\dot{\sigma})^2$—while maintaining second-order accuracy in the low-frequency range of interest \cite{chung1993time}.
\end{enumerate}
The following section details the derivation of the weak form and the linearization of the resulting nonlinear algebraic system.

\section{Finite element discretization and time integration}\label{sec:numerical}

To solve the nonlinear hyperbolic equation derived in \eqref{stress_wave_eq} numerically, we employ a semi-discrete approach. We first discretize the spatial domain using the Finite Element Method (FEM) and subsequently discretize the temporal domain using the Hilbert-Hughes-Taylor (HHT-$\alpha$) implicit integration scheme.

\subsection{Weak formulation}
Let $\Omega = (0, L)$ be the spatial domain. We denote by $H^1(\Omega)$ the standard Sobolev space containing functions with square-integrable derivatives. The space of test functions $\mathcal{V}$ consists of functions that satisfy the homogeneous form of the Dirichlet boundary conditions. Since we prescribed stress values at both $x=0$ and $x=L$, the test functions must vanish at these boundaries:
\begin{equation}
    \mathcal{V} = \{ w \in H^1(\Omega) \mid w(0) = 0 \text{ and } w(L) = 0 \}.
\end{equation}
The space of trial functions $\mathcal{S}_t$ depends on time due to the dynamic boundary loading. It consists of functions that satisfy the specific instantaneous boundary conditions given in the problem formulation:
\begin{equation}
    \mathcal{S}_t = \{ \sigma(\cdot, t) \in H^1(\Omega) \mid \sigma(0, t) = 0 \text{ and } \sigma(L, t) = A \sin(\omega t) \}.
\end{equation}
Note that $\sigma \in \mathcal{S}_t$ implies that $\sigma = v + \bar{\sigma}$, where $v \in \mathcal{V}$ and $\bar{\sigma}$ is a lift function satisfying the boundary data.  Multiplying Eq.~\eqref{stress_wave_eq} by an arbitrary test function $w \in \mathcal{V}$ and integrating over the domain yields:
\begin{equation}
    \int_{\Omega} w \rho \frac{\partial^2}{\partial t^2} [\mathfrak{f}(\sigma)] \, dx = \int_{\Omega} w \frac{\partial^2 \sigma}{\partial x^2} \, dx.
\end{equation}
Applying integration by parts to the spatial term on the right-hand side (and neglecting boundary terms for brevity, assuming homogeneous Neumann or Dirichlet conditions), we obtain the weak form: Find $\sigma \in \mathcal{S}$ such that for all $w \in \mathcal{V}$:
\begin{equation} \label{weak_form}
    \int_{\Omega} w \rho \left( \mathfrak{f}'(\sigma) \ddot{\sigma} + \mathfrak{f}''(\sigma) \dot{\sigma}^2 \right) \, dx + \int_{\Omega} \frac{\partial w}{\partial x} \frac{\partial \sigma}{\partial x} \, dx = 0.
\end{equation}
Note that we have utilized the expansion of the nonlinear inertial term derived in Eq.~\eqref{expanded_wave_eq}. This reveals that the weak form contains a state-dependent mass term and a convective-like velocity term.

\subsection{Existence and regularity of weak solutions}

The analysis of the existence of solutions for the variational problem \eqref{weak_form} presents significant analytical challenges compared to standard semi-linear wave equations. The weak form can be rewritten by identifying the effective mass density and the nonlinear source term:
\begin{equation}
    \langle \rho \mathfrak{f}'(\sigma) \ddot{\sigma}, w \rangle + \langle \nabla \sigma, \nabla w \rangle = - \langle \rho \mathfrak{f}''(\sigma) \dot{\sigma}^2, w \rangle \quad \forall w \in H^1_0(\Omega),
\end{equation}
where $\langle \cdot, \cdot \rangle$ denotes the $L^2$ inner product.

\subsubsection{Local existence via fixed-point arguments}
Local existence in time (for $t \in [0, T^*)$) can typically be established using a contraction mapping argument or the Galerkin method, provided the initial data is sufficiently smooth. The core requirement is the strict hyperbolicity condition:
\begin{equation}
    \mathfrak{f}'(\sigma) \geq c_0 > 0,
\end{equation}
for some constant $c_0$. Under this condition, the operator on the left-hand side retains the structure of a linear wave operator perturbed by state-dependent coefficients \cite{hughes1977well}. Standard theory for quasilinear hyperbolic systems implies that a unique classical solution exists for a short time interval, depending on the $H^s$-norm of the initial data ($\sigma_0, v_0$) \cite{kato1975cauchy,majda2012compressible}.

\subsubsection{Global existence and blow-up phenomena}
Global existence (for all $t > 0$) is not guaranteed and is, in fact, unlikely for arbitrary initial data due to the inherent structure of the nonlinearity.
\begin{itemize}
    \item \textbf{Shock formation:} Even with smooth initial conditions, quasilinear hyperbolic equations of this type are known to develop singularities (shocks) in finite time. This occurs because the wave characteristics propagate at speeds $c(\sigma)$ that depend on the solution amplitude. If characteristics converge, the gradients $\nabla \sigma$ and $\dot{\sigma}$ become unbounded, violating the regularity required for strong solutions.
    \item \textbf{The quadratic source term:} The term $\mathfrak{f}''(\sigma) \dot{\sigma}^2$ on the right-hand side acts as a nonlinear driving force. If $\mathfrak{f}''$ has a specific sign relative to the stress evolution, this term can lead to finite-time blow-up of the time derivative $\dot{\sigma}$, similar to the behavior observed in the Riccati ordinary differential equation.
\end{itemize}

Consequently, we seek \textit{weak} solutions in the distributional sense that satisfy an entropy condition to select the physically relevant solution after shock formation. For the numerical approximation presented in this work, the focus is on the pre-shock regime or the regularized solution obtained via the numerical dissipation inherent in the HHT-$\alpha$ time-integration scheme.
 
\subsection{Semi-discrete finite element system}
We approximate the stress field $\sigma(x,t)$ using standard finite element shape functions $N_I(x)$:
\begin{equation}
    \sigma(x,t) \approx \sigma^h(x,t) = \sum_{I=1}^{n_{dof}} N_I(x) \Sigma_I(t),
\end{equation}
where $\Sigma_I(t)$ are the time-dependent nodal stress values. Let $\boldsymbol{\Sigma}$, $\dot{\boldsymbol{\Sigma}}$, and $\ddot{\boldsymbol{\Sigma}}$ denote the global vectors for nodal stress, stress rate, and stress acceleration, respectively. Substituting the discretization into \eqref{weak_form} yields the nonlinear semi-discrete system of ordinary differential equations:
\begin{equation} \label{matrix_system}
    \mathbf{M}(\boldsymbol{\Sigma}) \ddot{\boldsymbol{\Sigma}} + \mathbf{F}_{nl}^{vel}(\boldsymbol{\Sigma}, \dot{\boldsymbol{\Sigma}}) + \mathbf{K} \boldsymbol{\Sigma} = \mathbf{0}.
\end{equation}
Here, the global matrices and vectors are defined as:
\begin{align}
    \mathbf{K}_{IJ} &= \int_{\Omega} \frac{\partial N_I}{\partial x} \frac{\partial N_J}{\partial x} \, dx \quad \text{(Constant stiffness matrix)}, \\
    \mathbf{M}_{IJ}(\boldsymbol{\Sigma}) &= \int_{\Omega} \rho \mathfrak{f}'(\sigma^h) N_I N_J \, dx \quad \text{(State-dependent mass matrix)}, \\
    \mathbf{F}_{nl, I}^{vel} &= \int_{\Omega} \rho \mathfrak{f}''(\sigma^h) (\dot{\sigma}^h)^2 N_I \, dx \quad \text{(Quadratic velocity force)}.
\end{align}
Unlike standard elastodynamics, where the nonlinearity resides in $\mathbf{K}$, here the stiffness matrix $\mathbf{K}$ is constant (linear), while the inertial forces (mass and velocity terms) are highly nonlinear.

\subsection{HHT-$\alpha$ time discretization}
We employ the HHT-$\alpha$ method to introduce controllable numerical dissipation for high-frequency noise while maintaining second-order accuracy. The method integrates the system from time step $t_n$ to $t_{n+1}$.

The Newmark kinematic updates are used to relate displacements, velocities, and accelerations:
\begin{align}
    \boldsymbol{\Sigma}_{n+1} &= \boldsymbol{\Sigma}_n + \Delta t \dot{\boldsymbol{\Sigma}}_n + \frac{\Delta t^2}{2} \left[ (1-2\beta) \ddot{\boldsymbol{\Sigma}}_n + 2\beta \ddot{\boldsymbol{\Sigma}}_{n+1} \right], \\
    \dot{\boldsymbol{\Sigma}}_{n+1} &= \dot{\boldsymbol{\Sigma}}_n + \Delta t \left[ (1-\gamma) \ddot{\boldsymbol{\Sigma}}_n + \gamma \ddot{\boldsymbol{\Sigma}}_{n+1} \right].
\end{align}
The discrete momentum balance is enforced using the HHT-$\alpha$ modification. The Newmark integration parameters are chosen as:
\begin{equation}
   \beta = \frac{(1-\alpha)^2}{4}, \quad \gamma = \frac{1}{2} - \alpha \quad \mbox{with} \quad \alpha \in \left[-\frac{1}{3}, 0 \right].
\end{equation}
Since the nonlinearity is primarily in the inertial terms, we apply the $\alpha$-shift to the linear stiffness term to ensure stability, while treating the nonlinear inertial forces fully implicitly at $t_{n+1}$ to satisfy the constitutive consistency:
\begin{equation} \label{residual}
    \mathbf{R}_{n+1} := \mathbf{F}^{inrt}(\boldsymbol{\Sigma}_{n+1}, \dot{\boldsymbol{\Sigma}}_{n+1}, \ddot{\boldsymbol{\Sigma}}_{n+1}) + (1-\alpha) \mathbf{K} \boldsymbol{\Sigma}_{n+1} + \alpha \mathbf{K} \boldsymbol{\Sigma}_n = \mathbf{0},
\end{equation}
where the total inertial force vector is $\mathbf{F}^{inrt} = \mathbf{M}(\boldsymbol{\Sigma}_{n+1})\ddot{\boldsymbol{\Sigma}}_{n+1} + \mathbf{F}_{nl}^{vel}(\boldsymbol{\Sigma}_{n+1}, \dot{\boldsymbol{\Sigma}}_{n+1})$.

\subsubsection{Linearization and solution strategy}
Since $\mathbf{R}_{n+1}$ is a nonlinear function of the primary unknown $\ddot{\boldsymbol{\Sigma}}_{n+1}$ (or equivalently $\boldsymbol{\Sigma}_{n+1}$ depending on implementation), we employ a Newton-Raphson iterative scheme. The linearized equation at iteration $k$ is:
\begin{equation}
    \mathbf{R}(\boldsymbol{\Sigma}^k_{n+1}) + \left[ \frac{\partial \mathbf{R}}{\partial \boldsymbol{\Sigma}} \frac{\partial \boldsymbol{\Sigma}}{\partial \boldsymbol{\Sigma}} + \frac{\partial \mathbf{R}}{\partial \dot{\boldsymbol{\Sigma}}} \frac{\partial \dot{\boldsymbol{\Sigma}}}{\partial \boldsymbol{\Sigma}} + \frac{\partial \mathbf{R}}{\partial \ddot{\boldsymbol{\Sigma}}} \frac{\partial \ddot{\boldsymbol{\Sigma}}}{\partial \boldsymbol{\Sigma}} \right] \delta \boldsymbol{\Sigma} = \mathbf{0}.
\end{equation}
Because the mass matrix depends on $\boldsymbol{\Sigma}$, the consistent tangent matrix $\mathbf{S}_{tan}$ becomes significantly more complex than in standard formulations:
\begin{equation}
    \mathbf{S}_{tan} = c_1 \frac{\partial \mathbf{F}^{inrt}}{\partial \ddot{\boldsymbol{\Sigma}}} + c_2 \frac{\partial \mathbf{F}^{inrt}}{\partial \dot{\boldsymbol{\Sigma}}} + c_3 \left( \frac{\partial \mathbf{F}^{inrt}}{\partial \boldsymbol{\Sigma}} + (1-\alpha)\mathbf{K} \right),
\end{equation}
where $c_1, c_2, c_3$ are Newmark coefficients determined by $\Delta t, \beta, \gamma$. The term $\frac{\partial \mathbf{F}^{inrt}}{\partial \boldsymbol{\Sigma}}$ involves the derivative of the mass matrix itself, requiring the third derivative $\mathfrak{f}'''(\sigma)$, highlighting the high degree of nonlinearity in this stress-based formulation.

\subsubsection{Final algebraic system and update procedure}
To advance the solution from time step $t_n$ to $t_{n+1}$, we solve the linearized system of equations at each Newton-Raphson iteration $k$. The final discrete algebraic system is given by:
\begin{equation} \label{linear_system}
    \mathbf{S}_{tan}(\boldsymbol{\Sigma}^k_{n+1}) \, \delta \boldsymbol{\Sigma}^k = - \mathbf{R}(\boldsymbol{\Sigma}^k_{n+1}),
\end{equation}
where $\delta \boldsymbol{\Sigma}^k$ is the increment of the nodal stress vector. The updated solution is obtained via $\boldsymbol{\Sigma}^{k+1}_{n+1} = \boldsymbol{\Sigma}^k_{n+1} + \delta \boldsymbol{\Sigma}^k$.

The consistent tangent stiffness matrix $\mathbf{S}_{tan}$ is derived by expanding the derivatives of the residual. Utilizing the Newmark parameters $\beta$ and $\gamma$, the sensitivity of the acceleration and velocity with respect to the primary variable $\boldsymbol{\Sigma}$ is given by:
\begin{equation}
    \frac{\partial \ddot{\boldsymbol{\Sigma}}}{\partial \boldsymbol{\Sigma}} = \frac{1}{\beta \Delta t^2} \mathbf{I}, \quad \frac{\partial \dot{\boldsymbol{\Sigma}}}{\partial \boldsymbol{\Sigma}} = \frac{\gamma}{\beta \Delta t} \mathbf{I}.
\end{equation}
Substituting these into the linearization yields the fully expanded form of the effective stiffness matrix:
\begin{equation} \label{tangent_matrix_expanded}
    \mathbf{S}_{tan} = \underbrace{\frac{1}{\beta \Delta t^2} \mathbf{M}(\boldsymbol{\Sigma})}_{\text{Inertial Stiffness}} 
    + \underbrace{\frac{\gamma}{\beta \Delta t} \mathbf{C}_{nl}(\boldsymbol{\Sigma}, \dot{\boldsymbol{\Sigma}})}_{\text{Velocity Stiffness}} 
    + \underbrace{(1-\alpha) \mathbf{K}}_{\text{Linear Stiffness}} 
    + \underbrace{\mathbf{K}_{geo}(\boldsymbol{\Sigma}, \ddot{\boldsymbol{\Sigma}})}_{\text{Nonlinear Geometric Stiffness}},
\end{equation}
where the nonlinear contributions are defined as:
\begin{align}
    \mathbf{C}_{nl} &= \frac{\partial \mathbf{F}^{inrt}}{\partial \dot{\boldsymbol{\Sigma}}} = \int_{\Omega} 2 \rho \mathfrak{f}''(\sigma) \dot{\sigma} N_I N_J \, dx, \\
    \mathbf{K}_{geo} &= \frac{\partial \mathbf{M}}{\partial \boldsymbol{\Sigma}} \ddot{\boldsymbol{\Sigma}} = \int_{\Omega} \rho \mathfrak{f}'''(\sigma) \ddot{\sigma} N_I N_J \, dx.
\end{align}
Equation \eqref{tangent_matrix_expanded} highlights the complexity of the formulation: unlike standard linear elastodynamics where $\mathbf{S}_{tan}$ is constant (or state-dependent only via material stiffness), here the "mass" contribution varies with every iteration, and the acceleration term induces a higher-order "geometric" stiffness $\mathbf{K}_{geo}$ dependent on the third derivative of the constitutive function.

The iterative procedure continues until the norm of the residual $\|\mathbf{R}_{n+1}\|$ falls below a user-defined tolerance. Upon convergence, the kinematic variables $\dot{\boldsymbol{\Sigma}}$ and $\ddot{\boldsymbol{\Sigma}}$ are updated consistent with the Newmark relations.

\section{Numerical implementation and results}\label{sec:results}

The nonlinear stress wave equation derived in the previous sections was implemented using the open-source finite element library \texttt{deal.II} \cite{arndt2021deal}. This library was selected for its robust support of $hp$-adaptive finite element methods (FEM) and dimension-independent programming. 

\subsection{Computational framework and adaptivity}
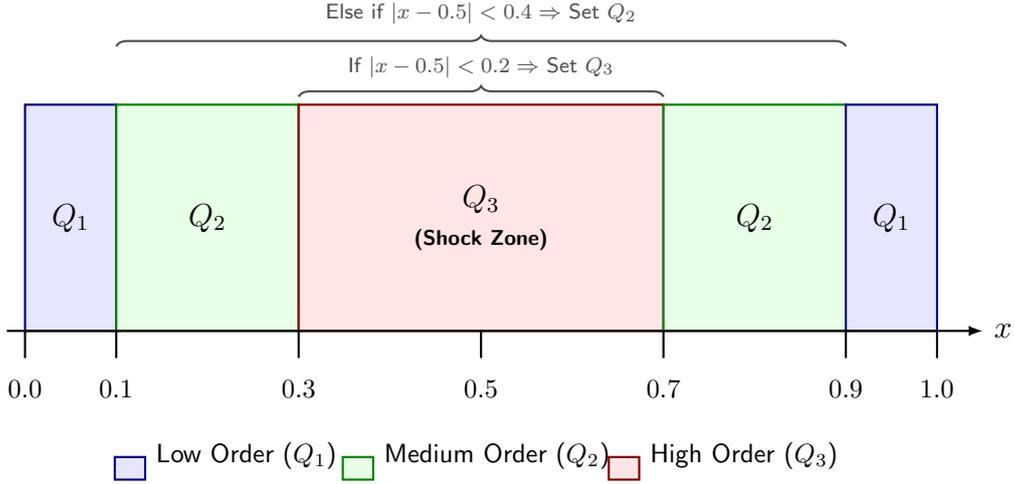
\begin{figure}[H]
    \centering
    
    \begin{tikzpicture}[
        scale=12, 
        styleQ1/.style={fill=blue!10, draw=blue!50!black, line width=0.8pt},
        styleQ2/.style={fill=green!10, draw=green!50!black, line width=0.8pt},
        styleQ3/.style={fill=red!10, draw=red!50!black, line width=0.8pt},
        labelStyle/.style={font=\sffamily\bfseries\normalsize, align=center},
        logicFont/.style={font=\sffamily\scriptsize, align=center, fill=white, inner sep=2pt}
    ]
        
        \draw[styleQ1] (0.0, 0) rectangle (0.1, 0.25) node[midway, labelStyle] {$Q_1$};
        
        \draw[styleQ2] (0.1, 0) rectangle (0.3, 0.25) node[midway, labelStyle] {$Q_2$};
        
        \draw[styleQ3] (0.3, 0) rectangle (0.7, 0.25) node[midway, labelStyle] {$Q_3$\\ \scriptsize (Shock Zone)};
        
        \draw[styleQ2] (0.7, 0) rectangle (0.9, 0.25) node[midway, labelStyle] {$Q_2$};
        
        \draw[styleQ1] (0.9, 0) rectangle (1.0, 0.25) node[midway, labelStyle] {$Q_1$};

        \draw[->, thick, >=latex] (-0.02, 0) -- (1.05, 0) node[right, font=\bfseries] {$x$};
        
        \foreach \x in {0.0, 0.1, 0.3, 0.5, 0.7, 0.9, 1.0} {
            \draw[thick] (\x, 0) -- (\x, -0.03) node[below=4pt, font=\sffamily\footnotesize] {$\x$};
        }

        \draw[decorate, decoration={brace, amplitude=4pt, raise=3pt}, thick, black!70] 
            (0.3, 0.25) -- (0.7, 0.25) 
            node[midway, above=8pt, logicFont] {If $|x - 0.5| < 0.2 \Rightarrow$ Set $Q_3$};

        \draw[decorate, decoration={brace, amplitude=4pt, raise=22pt}, thick, black!70] 
            (0.1, 0.25) -- (0.9, 0.25) 
            node[midway, above=28pt, logicFont] {Else if $|x - 0.5| < 0.4 \Rightarrow$ Set $Q_2$};

        \node[below=1.2cm] at (0.5, 0) {
            \begin{tikzpicture}[scale=1]
                \node[styleQ1, minimum width=0.4cm, minimum height=0.3cm] at (0,0) {};
                \node[right, font=\sffamily\footnotesize] at (0.2,0) {Low Order ($Q_1$)};
                
                \node[styleQ2, minimum width=0.4cm, minimum height=0.3cm] at (3.0,0) {};
                \node[right, font=\sffamily\footnotesize] at (3.2,0) {Medium Order ($Q_2$)};
                
                \node[styleQ3, minimum width=0.4cm, minimum height=0.3cm] at (6.5,0) {};
                \node[right, font=\sffamily\footnotesize] at (6.7,0) {High Order ($Q_3$)};
            \end{tikzpicture}
        };

    \end{tikzpicture}
    
    \caption{Adaptive polynomial degree distribution ($h$-$p$ refinement strategy). The domain is divided based on distance from the center ($L/2$). High-order $Q_3$ elements are localized in the center to resolve shock formation, transitioning to $Q_1$ at the boundaries.}
    \label{fig:poly_logic}
\end{figure}

A primary challenge in solving the stress-based wave equation \eqref{stress_wave_eq} is the potential formation of shock-like gradients in the stress variable, even from smooth initial data, due to the material nonlinearity. To address this, our implementation utilizes a spatially adaptive mesh strategy. As illustrated in Figure~\ref{fig:poly_logic}, we employ an $hp$-refinement approach where the polynomial degree of the finite element basis functions ($p$) varies across the domain. High-order cubic elements ($Q_3$) are localized in regions anticipating high gradients (such as the shock zone in the domain center), while lower-order linear elements ($Q_1$) are utilized in the far field to reduce computational cost.

\subsection{Nonlinear solution and kinematic reconstruction}
The global nonlinear system is solved using a Newton-Raphson scheme within each HHT-$\alpha$ time step. A unique feature of this pure stress formulation is the absence of displacement degrees of freedom in the primary system solution. Consequently, physical kinematic quantities—displacement $u$ and particle velocity $v$—must be reconstructed via spatial integration of the constitutive relations during the post-processing phase. 

Algorithm~\ref{alg:fem_solver} details the complete solution procedure. Note specifically lines 23-31, where the strain $\varepsilon(\sigma)$ and strain-rate $\dot{\varepsilon}(\sigma, \dot{\sigma})$ are integrated to recover the motion of the body. This reconstruction serves as a consistency check for the formulation.

\begin{algorithm}[H]
\caption{Nonlinear Finite Element Solver with Post-Processing}
\label{alg:fem_solver}
\begin{algorithmic}[1]
\Require Mesh $\mathcal{T}$, Parameters $\beta, \alpha, \Delta t, T_{final}$, Tol
\State Initialize $\sigma_0 = 0, \dot{\sigma}_0 = 0, \ddot{\sigma}_0 = 0$
\State $t \gets 0$

\While{$t < T_{final}$} \Comment{Time Stepping Loop}
    \State $t \gets t + \Delta t$
    \State $\ddot{\sigma}_{n+1}^{(0)} \gets \ddot{\sigma}_n$ \Comment{Initial guess for acceleration}
    \State $k \gets 0$
    
    \Repeat \Comment{Newton-Raphson Iteration}
        \State \textbf{1. Newmark Predictor:}
        \State $\sigma_{n+1}^{(k)} \gets \sigma_n + \Delta t \dot{\sigma}_n + \Delta t^2 (0.5 - \beta_{nm}) \ddot{\sigma}_n + \Delta t^2 \beta_{nm} \ddot{\sigma}_{n+1}^{(k)}$
        \State $\dot{\sigma}_{n+1}^{(k)} \gets \dot{\sigma}_n + \Delta t (1 - \gamma_{nm}) \ddot{\sigma}_n + \Delta t \gamma_{nm} \ddot{\sigma}_{n+1}^{(k)}$

        \State \textbf{2. Assembly (Weak Form):}
        \State Compute Tangent Stiffness $K$ and Residual $R$:
        \State $C_t(\sigma) \gets \frac{d\varepsilon}{d\sigma} = (1 + (\beta|\sigma|)^\alpha)^{-(1 + 1/\alpha)}$
        \State $R \gets \int_{\Omega} \left[ \rho \ddot{\varepsilon}(\sigma) \phi + \nabla \sigma \cdot \nabla \phi \right] d\Omega$
        
        \State \textbf{3. Linear Solve:}
        \State Solve $K \Delta \ddot{\sigma} = -R$
        \State Update: $\ddot{\sigma}_{n+1}^{(k+1)} \gets \ddot{\sigma}_{n+1}^{(k)} + \Delta \ddot{\sigma}$
        \State $k \gets k + 1$
    \Until{$||R|| < \text{Tol}$ or $k > k_{max}$}

    \State Update history: $\sigma_n \gets \sigma_{n+1}, \dot{\sigma}_n \gets \dot{\sigma}_{n+1}, \ddot{\sigma}_n \gets \ddot{\sigma}_{n+1}$

    \If{Output Step} \Comment{Post-Processing Routine}
        \State Define sampling points $x_i$ for $i=0 \dots M$
        \State Evaluate $\sigma(x_i)$ and $\dot{\sigma}(x_i)$ from FE solution
        \State Initialize $u(0) = 0, v(0) = 0$
        
        \For{$i = 1$ to $M$} \Comment{Spatial Integration}
            \State \textbf{Compute Strain:}
            \State $\varepsilon_i \gets \frac{\sigma_i}{(1 + (\beta|\sigma_i|)^\alpha)^{1/\alpha}}$
            \State \textbf{Integrate Displacement ($u = \int \varepsilon dx$):}
            \State $u_i \gets u_{i-1} + \frac{1}{2}(\varepsilon_i + \varepsilon_{i-1}) \Delta x$
            \State \textbf{Integrate Particle Velocity ($v = \int \dot{\varepsilon} dx$):}
            \State $\dot{\varepsilon}_i \gets C_t(\sigma_i) \cdot \dot{\sigma}_i$
            \State $v_i \gets v_{i-1} + \frac{1}{2}(\dot{\varepsilon}_i + \dot{\varepsilon}_{i-1}) \Delta x$
            \State \textbf{Compute Local Wave Speed:}
            \State $c_i \gets \sqrt{\frac{1}{\rho \cdot C_t(\sigma_i)}}$
        \EndFor
        \State Write $\{x_i, \sigma_i, u_i, v_i, \varepsilon_i, c_i\}$ to file
    \EndIf
\EndWhile
\end{algorithmic}
\end{algorithm} 

\subsection{Verification and Convergence Analysis}
Validation of the numerical framework is conducted using the Method of Manufactured Solutions (MMS) \cite{roache2002code}. This technique allows for the rigorous assessment of the code's ability to recover theoretical convergence rates for a known smooth solution.

To verify the order of accuracy, we manufacture an exact solution for the stress field:
\begin{equation}
    \sigma_{\text{exact}}(x,t) = \sin(\pi x) \sin(t), \quad x \in [0, 1]
\end{equation}
This solution satisfies homogeneous Dirichlet boundary conditions $\sigma(0,t) = \sigma(1,t) = 0$. The model parameters used for the verification test are:
\begin{itemize}
    \item Density: $\rho = 1.0$
    \item Material parameters: $b = 1.0$, $a = 2.0$ (Smooth constitutive law)
    \item HHT Damping: $\alpha = -0.05$
\end{itemize}

\subsubsection{Spatial Convergence (Order of Accuracy: 2)}
The spatial discretization error is isolated by fixing a sufficiently small time step $\Delta t = 10^{-5}$ to eliminate temporal integration errors. We vary the mesh refinement level ($h$) using linear Lagrange elements ($Q_1$). Theoretical finite element analysis predicts an asymptotic convergence rate of $O(h^{p+1})$ for the $L_2$ norm, which corresponds to a rate of 2.0 for $p=1$.

The results in Table \ref{tab:spatial} confirm this theoretical prediction. As the mesh is refined from 16 to 128 cells, the computed rate stabilizes at exactly 2.00, demonstrating that the spatial weak form and tangent matrix assembly are correctly implemented.

\begin{table}[H]
    \centering
    \caption{Spatial convergence rates using $Q_1$ elements ($\Delta t = 10^{-5}$).}
    \label{tab:spatial}
    \begin{tabular}{c c c c}
        \toprule
        \textbf{Refinement Level} & \textbf{DoFs} & \textbf{$L_2$ Error} & \textbf{Rate} \\
        \midrule
        4 (16 cells)  & 17  & \num{1.246e-04} & --   \\
        5 (32 cells)  & 33  & \num{3.117e-05} & 2.00 \\
        6 (64 cells)  & 65  & \num{7.793e-06} & 2.00 \\
        7 (128 cells) & 129 & \num{1.948e-06} & 2.00 \\
        \bottomrule
    \end{tabular}
\end{table}

\subsubsection{Temporal Convergence (Order of Accuracy: 2)}
Temporal accuracy is assessed by fixing a high spatial refinement and using high-order Cubic elements ($Q_3$) to render spatial errors negligible compared to time integration errors. The HHT-$\alpha$ method is designed to be second-order accurate, $O(\Delta t^2)$, even for nonlinear problems.

Table \ref{tab:temporal} presents the $L_2$ error norm as the time step $\Delta t$ is halved. We observe an asymptotic convergence toward the theoretical rate of 2.0. The rate progression (1.40 $\to$ 1.84 $\to$ 1.92) is characteristic of nonlinear inertial problems solved with HHT-$\alpha$. At larger time steps, the linearization error of the highly nonlinear mass term $\mathbf{M}(\boldsymbol{\Sigma})$ dominates, slightly degrading the observed rate. However, as $\Delta t$ decreases, the second-order truncation error of the time integrator becomes dominant, and the rate approaches the optimal value of 2.0.

\begin{table}[H]
    \centering
    \caption{Temporal convergence rates using HHT-$\alpha$ ($\alpha_{\text{HHT}}=-0.05$).}
    \label{tab:temporal}
    \begin{tabular}{c c c}
        \toprule
        \textbf{Time Step ($\Delta t$)} & \textbf{$L_2$ Error} & \textbf{Rate} \\
        \midrule
        \num{8.000e-03} & \num{2.394e-06} & --   \\
        \num{4.000e-03} & \num{1.475e-06} & 1.40 \\
        \num{2.000e-03} & \num{7.795e-07} & 1.84 \\
        \num{1.000e-03} & \num{4.003e-07} & 1.92 \\
        \bottomrule
    \end{tabular}
\end{table}

\subsection{Linear case with $b=0$}

\begin{figure}[H]
    \centering
    \begin{subfigure}[b]{0.48\textwidth}
        \centering
        \includegraphics[width=\textwidth, height=9cm]{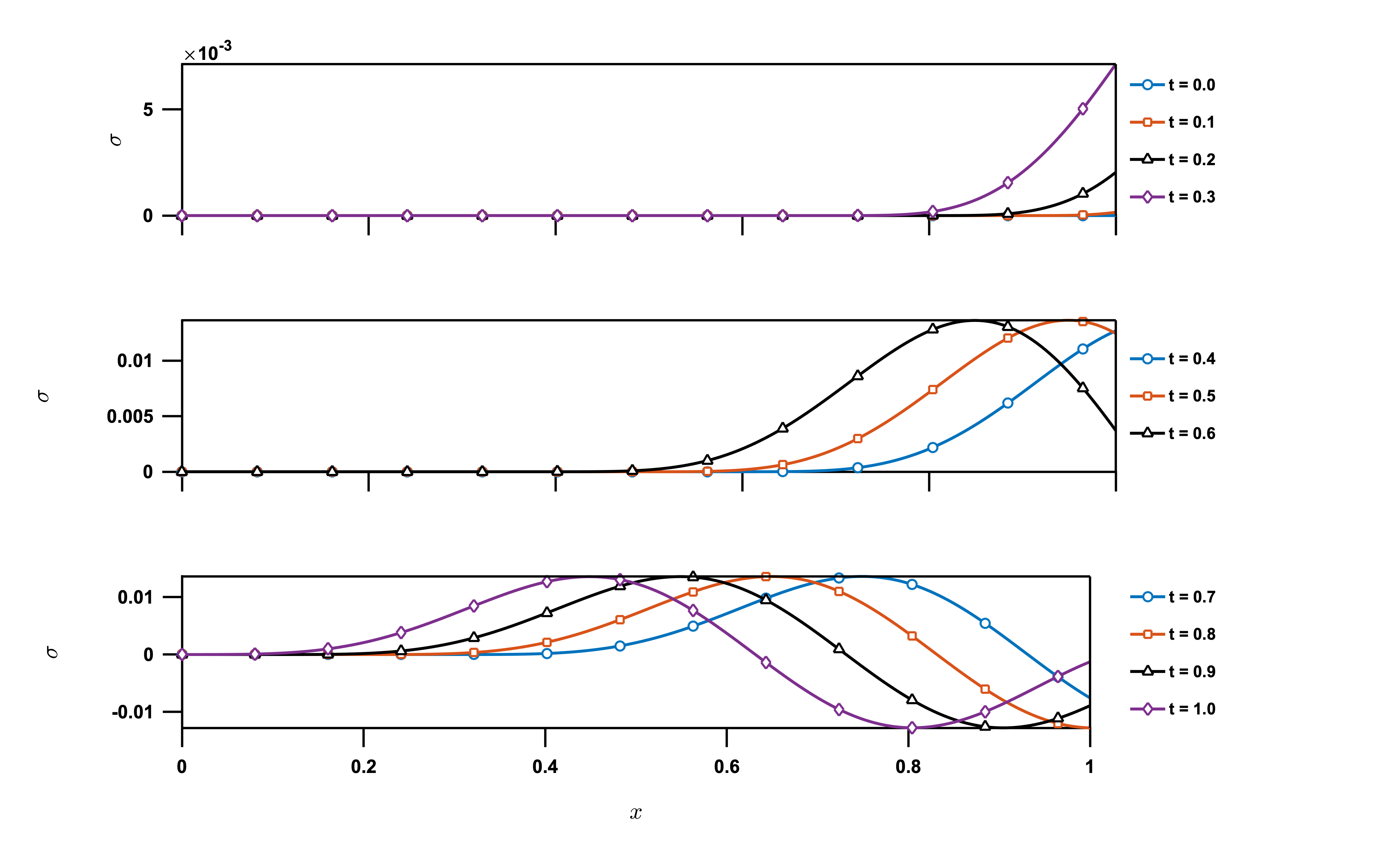}
        \caption{$\sigma$}
        \label{fig:sigma_lin}
    \end{subfigure}
    \hfill 
    \begin{subfigure}[b]{0.48\textwidth}
        \centering
        \includegraphics[width=\textwidth, height=9cm]{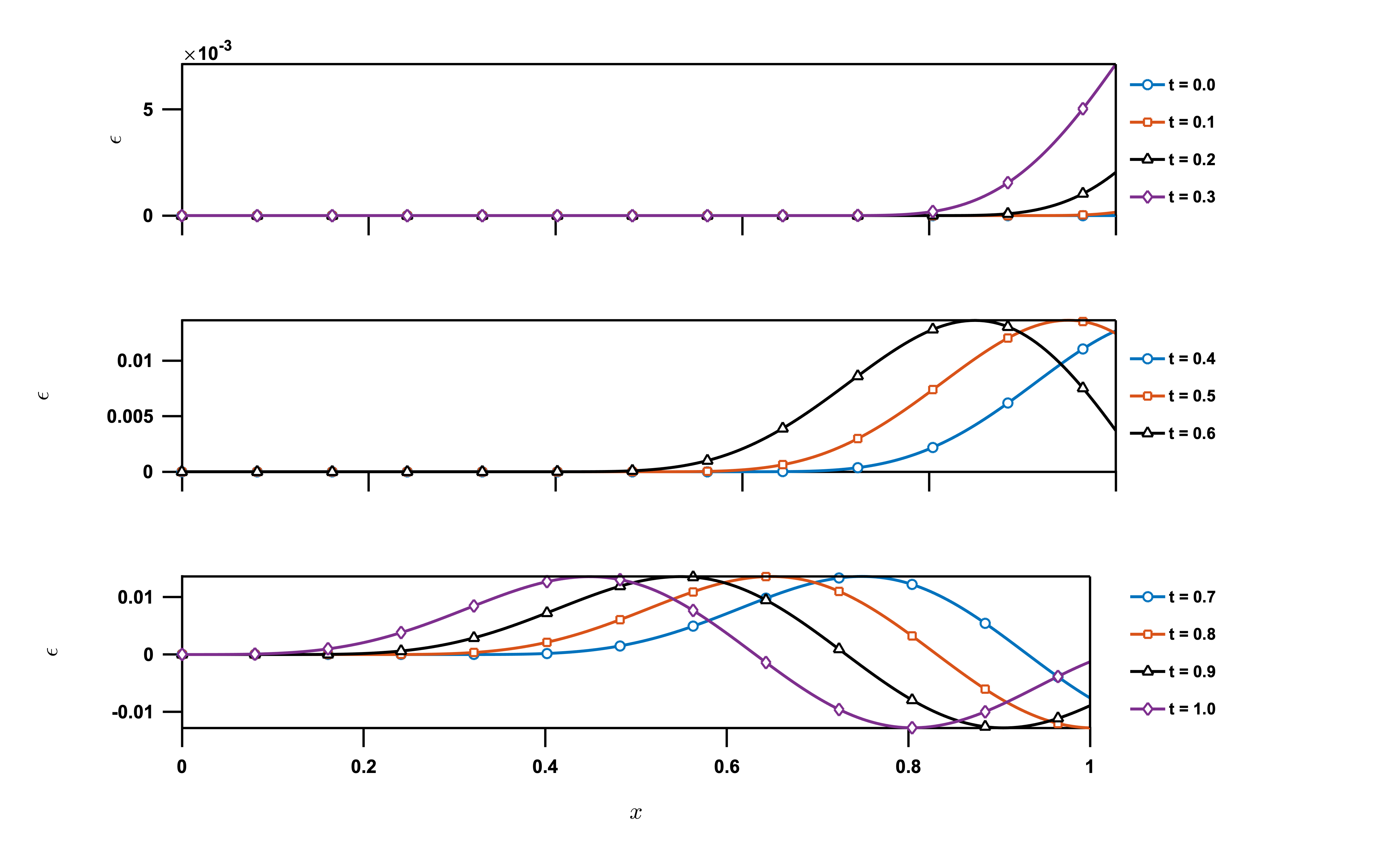}
        \caption{$\epsilon$}
        \label{fig:strain_lin}
    \end{subfigure}    
    \caption{The stress $\sigma$ and strain $\epsilon$ for the linear case with $b=0.0$ and $a=1.5$. Note the identical wave profiles characteristic of linear elasticity.}
    \label{fig:stress_strain_linear}
\end{figure}

\begin{figure}[H]
    \centering
    \begin{subfigure}[b]{0.48\textwidth}
        \centering
        \includegraphics[width=\textwidth, height=9cm]{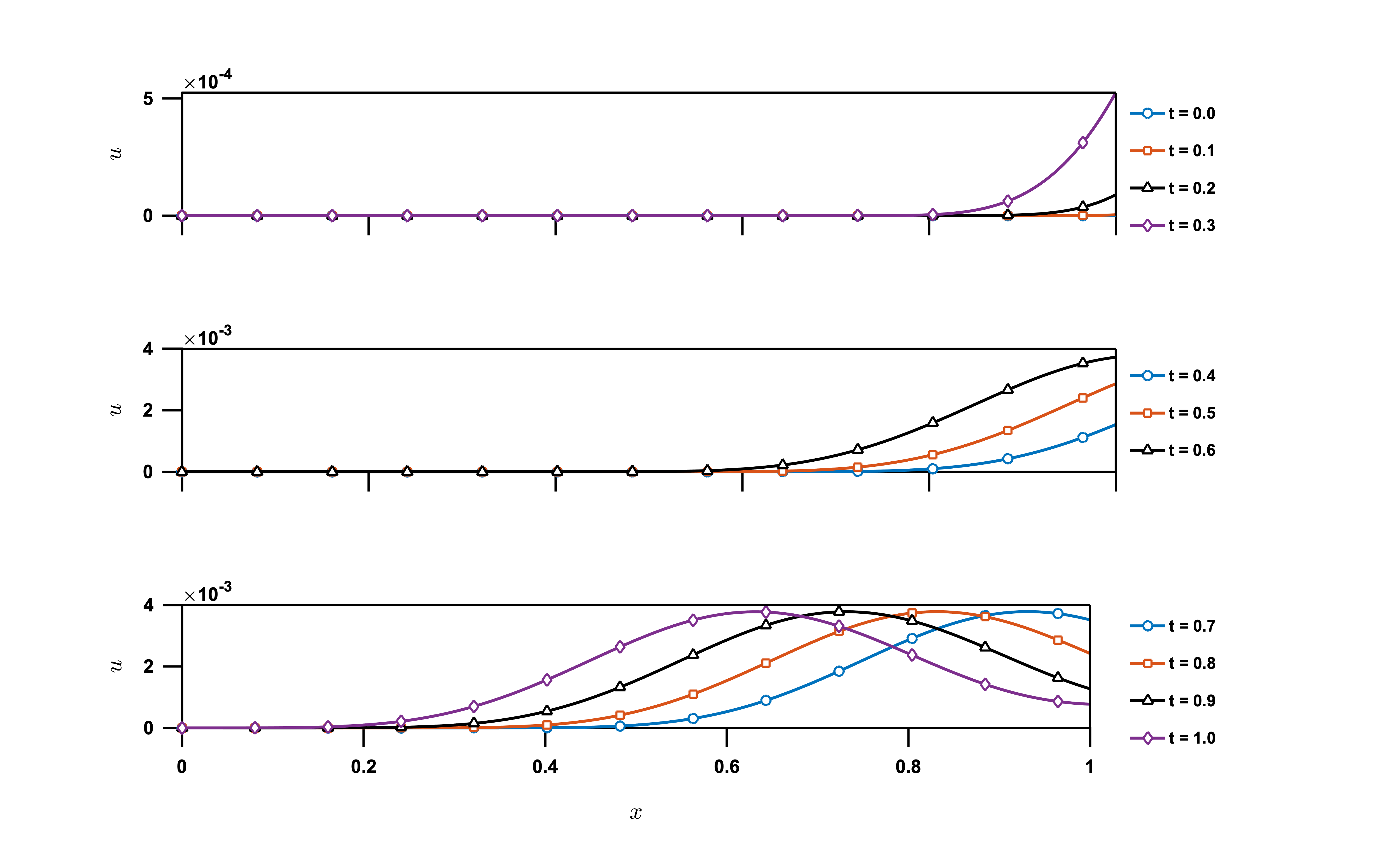}
        \caption{$u$}
        \label{fig:disp_lin}
    \end{subfigure}
    \hfill 
    \begin{subfigure}[b]{0.48\textwidth}
        \centering
        \includegraphics[width=\textwidth, height=9cm]{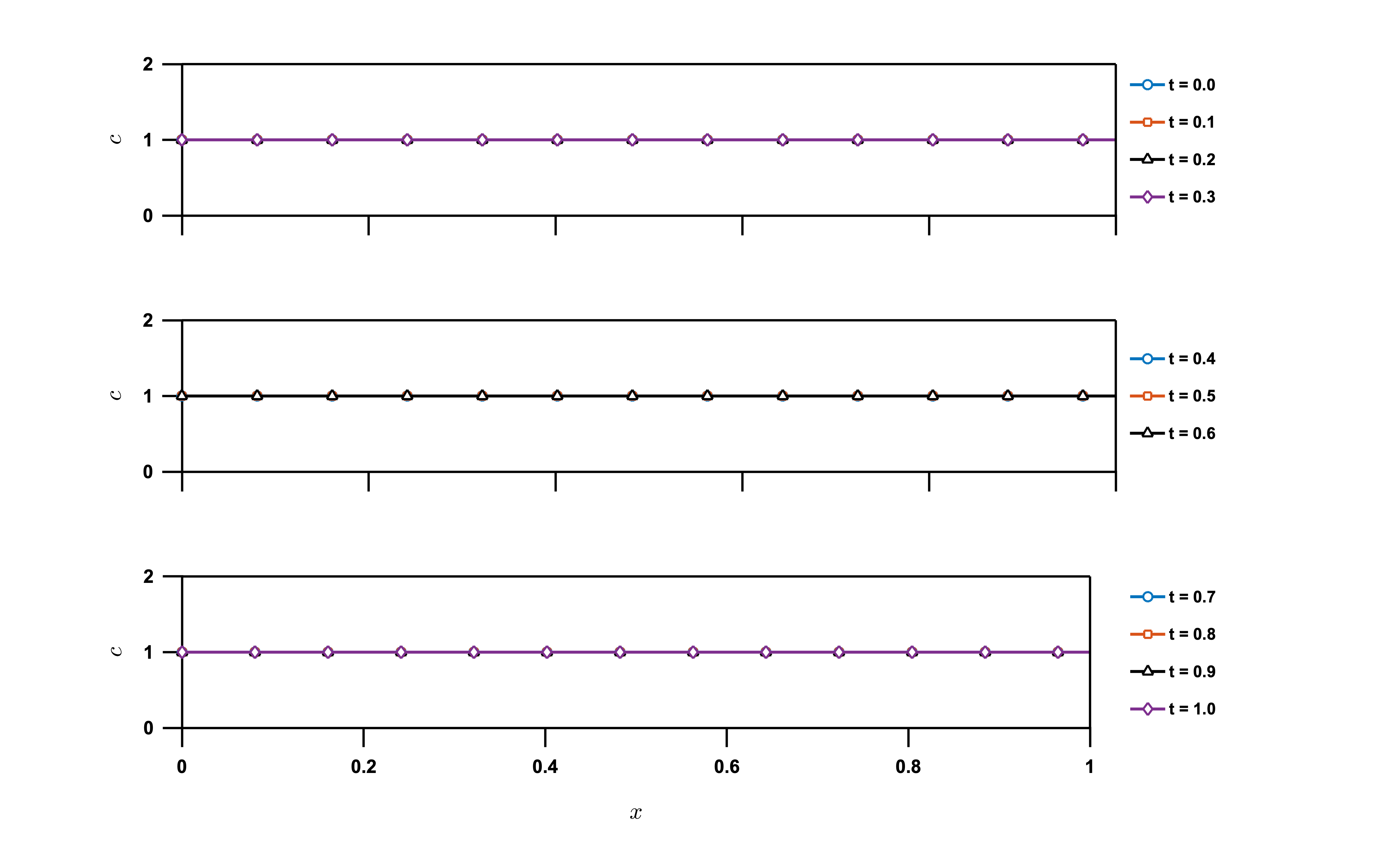}
        \caption{$c$}
        \label{fig:speed_lin}
    \end{subfigure}   
    \caption{The displacement $u$ and wave speed $c$ for the linear case with $b=0.0$ and $a=1.5$. The constant wave speed $c=1$ confirms the linearity of the medium.}
    \label{fig:disp_speed_linear}
\end{figure}

This subsection examines the linear elastic limit defined by parameters $b = 0.0$ and $a = 1.5$. As shown in Fig. \ref{fig:speed_lin}, the wave speed $c$ remains invariant at unity ($c=1$) throughout the spatiotemporal domain, confirming that the material stiffness is independent of the deformation state. Consequently, the displacement $u$, stress $\sigma$, and strain $\epsilon$ propagate as stable, non-dispersive waves. The pulse, initiated at the boundary $x=1$, travels leftward into the domain while preserving its smooth topology. Unlike the nonlinear cases discussed later, the stress and strain profiles (Fig. \ref{fig:stress_strain_linear}) remain identical in shape—differing only by the elastic modulus—and exhibit no gradient steepening or shock formation.

\subsection{Nonlinear case with $b=1.0$ and $a=1.5$}
\begin{figure}[H]
    \centering
    \begin{subfigure}[b]{0.48\textwidth}
        \centering
        \includegraphics[width=\textwidth, height=9cm]{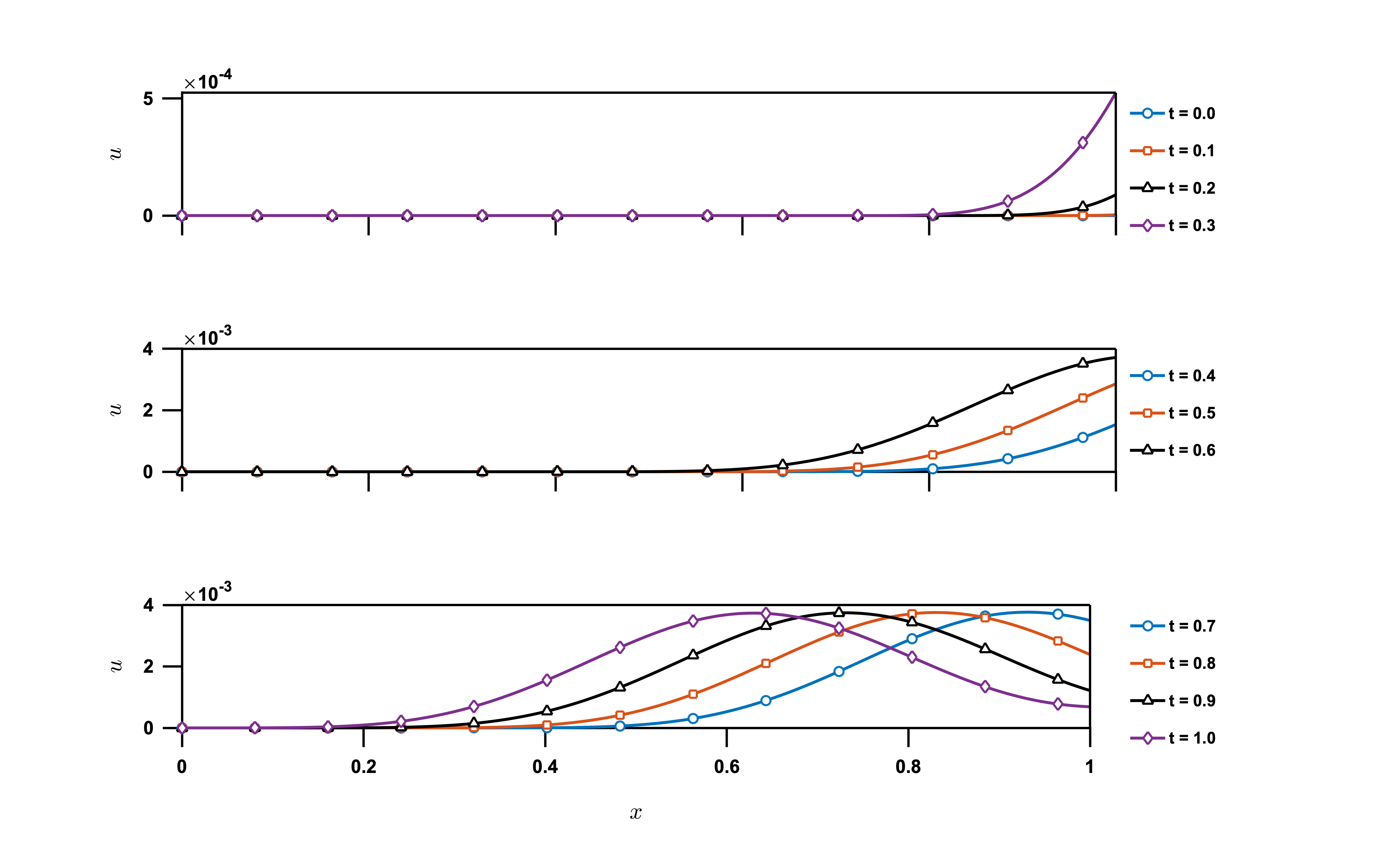}
        \caption{$u$}
        \label{fig:disp_beta1}
    \end{subfigure}
    \hfill 
    \begin{subfigure}[b]{0.48\textwidth}
        \centering
        \includegraphics[width=\textwidth, height=9cm]{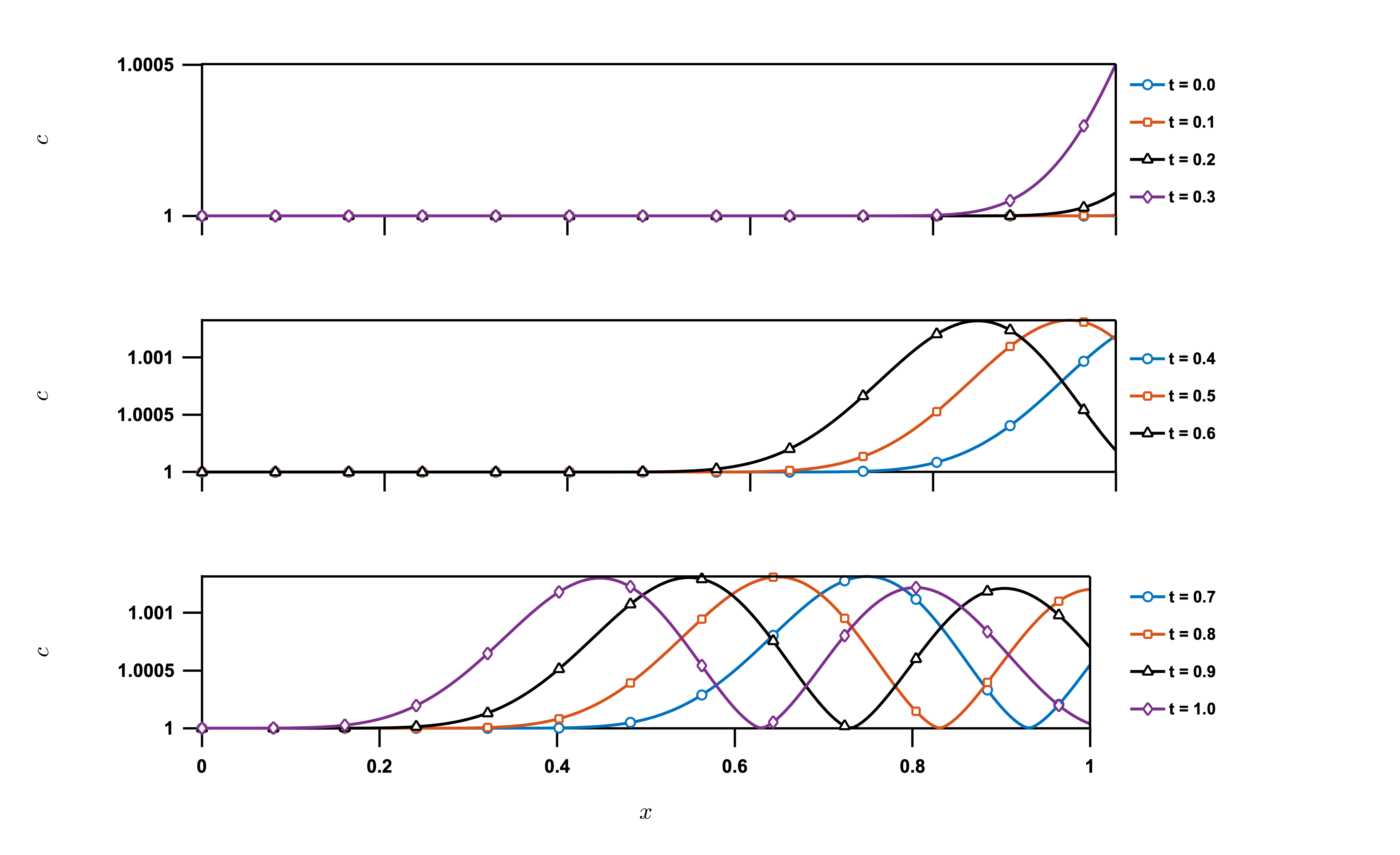}
        \caption{$c$}
        \label{fig:speed_beta1}
    \end{subfigure}   
    \caption{The displacement $u$ and wave speed $c$ for the nonlinear case with $b=1.0$ and $a=1.5$. Note the slight variation in wave speed $c$ compared to the constant value in the linear case.}
    \label{fig:disp_speed_beta1}
\end{figure}

In the weakly nonlinear regime defined by $b = 1.0$ and $a = 1.5$, the wave propagation dynamics diverge from the linear baseline.. In contrast to the linear regime ($b=0$), the wave speed $c$ is no longer constant. As shown in Fig. \ref{fig:speed_beta1}, $c$ exhibits small temporal and spatial variations, fluctuating slightly above unity. This variation indicates that the material stiffness is now dependent on the local deformation state. Consequently, the mechanical fields—stress $\sigma$, strain $\epsilon$, and displacement $u$—begin to show subtle deviations from the purely linear solution. Although the nonlinearity is mild at $b=1.0$, the dependence of wave speed on strain suggests the onset of dispersive effects.

\begin{figure}[H]
    \centering
    \begin{subfigure}[b]{0.48\textwidth}
        \centering
        \includegraphics[width=\textwidth, height=9cm]{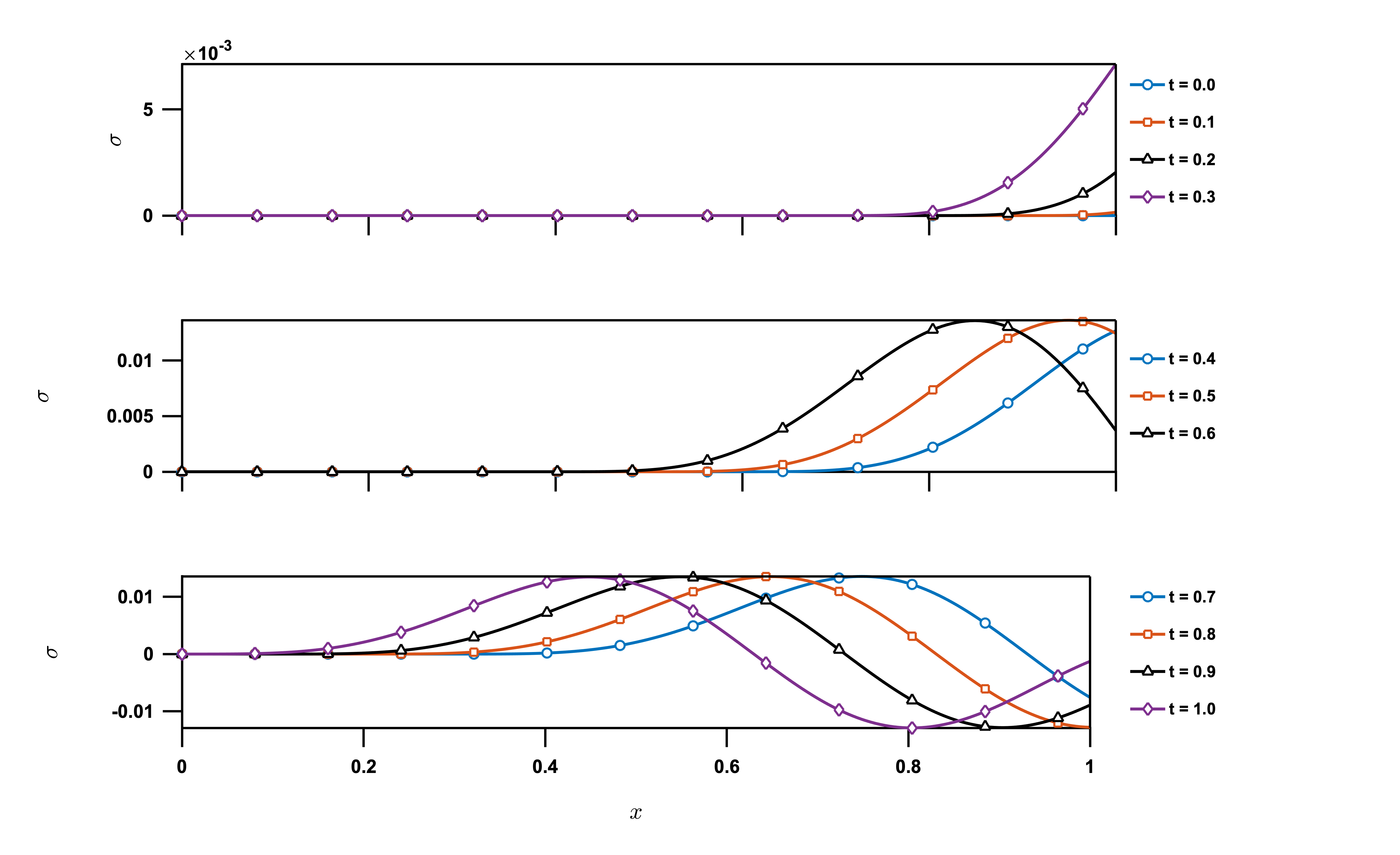}
        \caption{$\sigma$}
        \label{fig:sigma_beta1}
    \end{subfigure}
    \hfill 
    \begin{subfigure}[b]{0.48\textwidth}
        \centering
        \includegraphics[width=\textwidth, height=9cm]{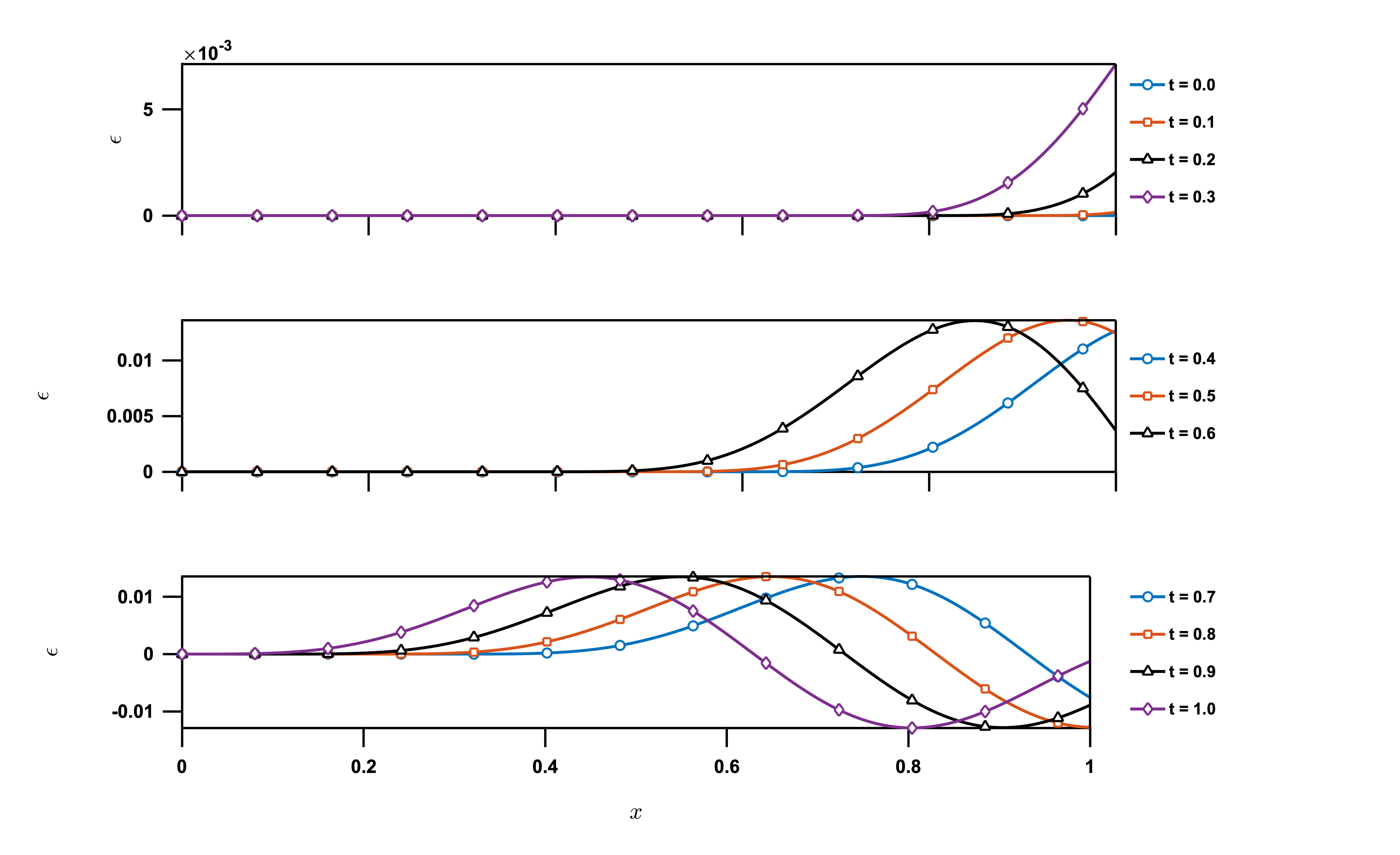}
        \caption{$\epsilon$}
        \label{fig:strain_beta1}
    \end{subfigure}    
    \caption{The stress $\sigma$ and strain $\epsilon$ for the weakly nonlinear case with $b=1.0$ and $a=1.5$.}
    \label{fig:stress_strain_beta1}
\end{figure}

\subsection{Nonlinear case with $b=5.0$}

\begin{figure}[H]
    \centering

    \begin{subfigure}[b]{0.48\textwidth}
        \centering
        \includegraphics[width=\textwidth, height=9cm]{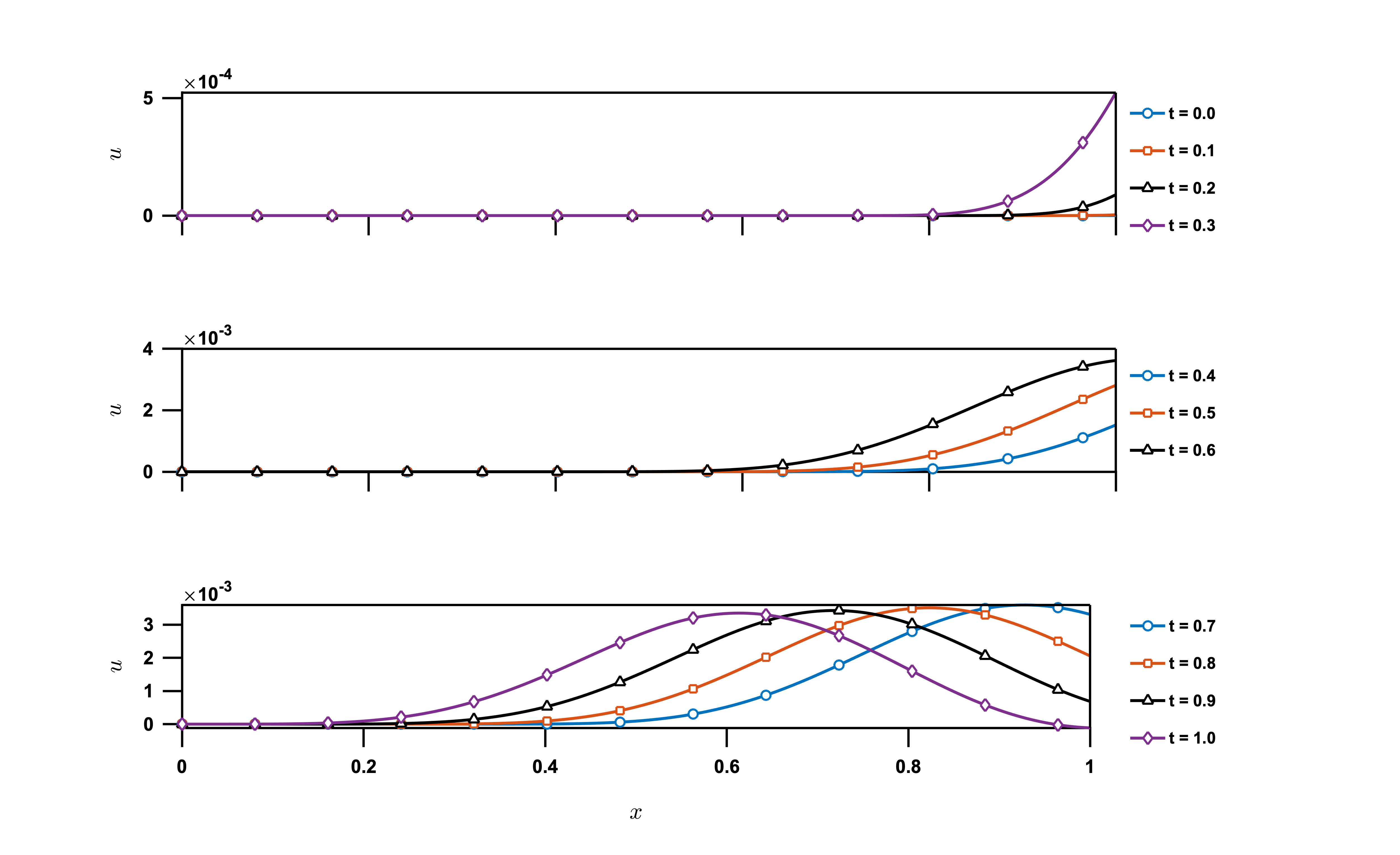}
        \caption{$u$}
        \label{fig:disp_beta5}
    \end{subfigure}
    \hfill 
    \begin{subfigure}[b]{0.48\textwidth}
        \centering
        \includegraphics[width=\textwidth, height=9cm]{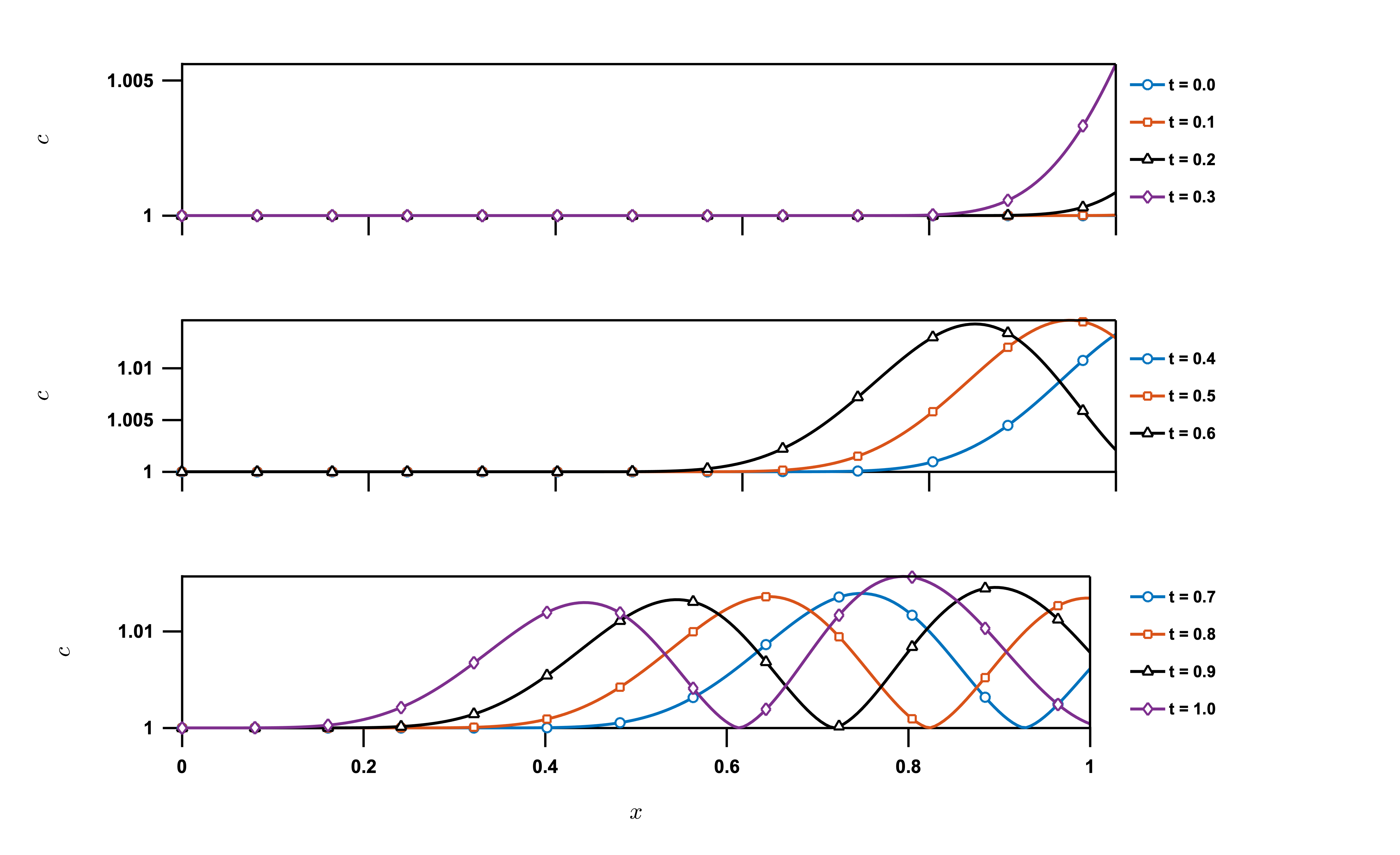}
        \caption{$c$}
        \label{fig:speed_beta5}
    \end{subfigure}
    \caption{The displacement $u$ and wave speed $c$ for the nonlinear case with $\beta=5.0$ and $\alpha=1.5$. Note the significant increase in wave speed variation compared to lower $\beta$ values.}
    \label{fig:disp_speed_beta5}
\end{figure}

We now consider the moderately nonlinear regime characterized by $b = 5.0$ and $\alpha = 1.5$. In this case, the influence of the constitutive nonlinearity becomes significantly more pronounced compared to the weakly nonlinear scenario. The wave speed $c$, illustrated in Fig. \ref{fig:speed_beta5}, exhibits substantial spatial and temporal fluctuations, reaching peak values of approximately $c \approx 1.015$. This magnitude of variation—roughly an order of magnitude larger than in the $b=1.0$ case—demonstrates a strong coupling between the local deformation state and the material stiffness. As a direct consequence of this stress-dependent wave velocity, the propagating pulses for stress $\sigma$ and strain $\epsilon$ undergo visible distortion as they traverse the domain. Unlike the linear case where the pulse shape is preserved, the higher amplitude regions of the wave now travel faster than the lower amplitude tails. This differential velocity leads to a gradual steepening of the wavefront, a precursor to shock formation, which is a hallmark behavior of nonlinear hyperbolic conservation laws.

\begin{figure}[H]
    \centering

    \begin{subfigure}[b]{0.48\textwidth}
        \centering
        \includegraphics[width=\textwidth, height=9cm]{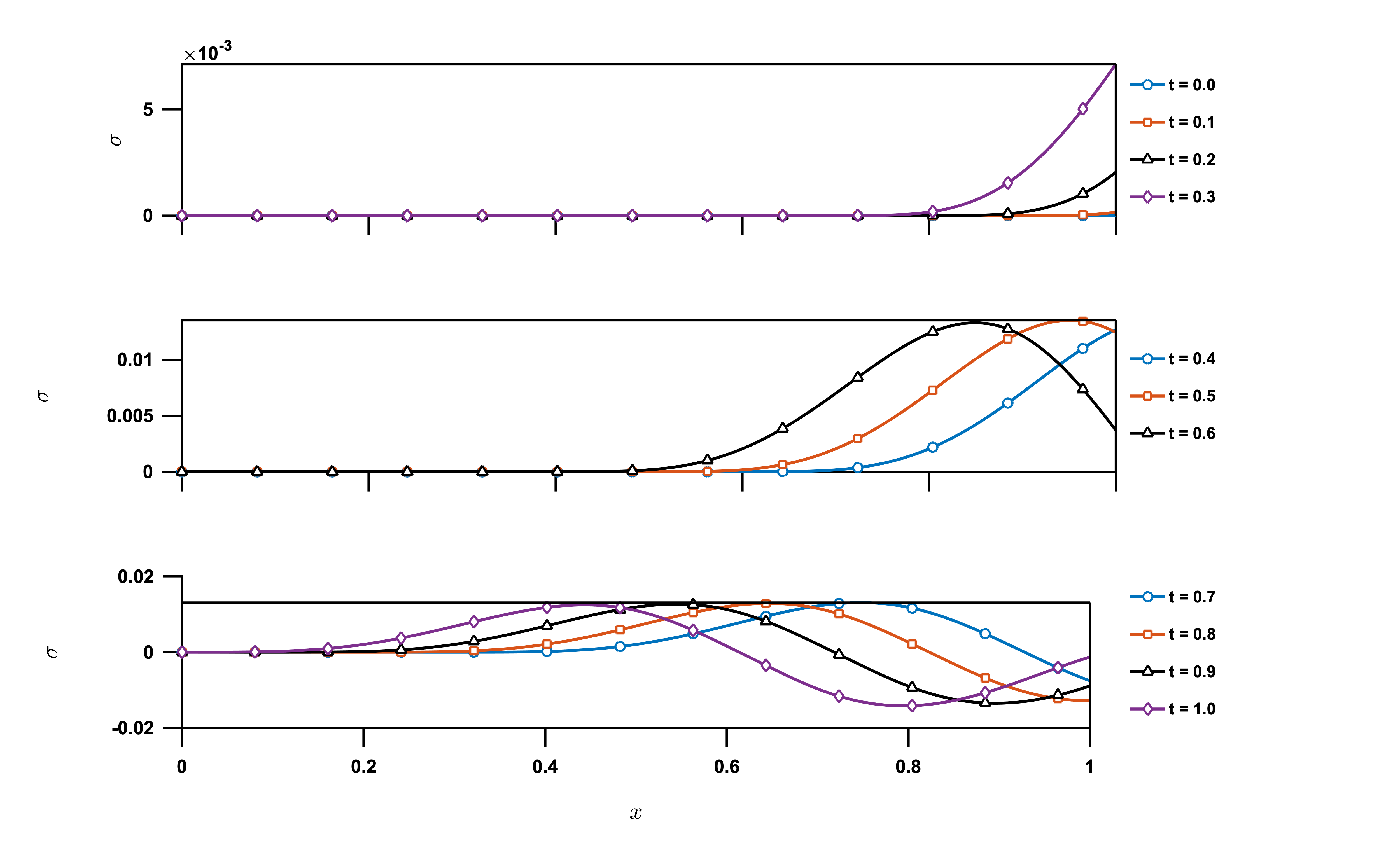}
        \caption{$\sigma$}
        \label{fig:sigma_beta5}
    \end{subfigure}
    \hfill 
    \begin{subfigure}[b]{0.48\textwidth}
        \centering
        \includegraphics[width=\textwidth, height=9cm]{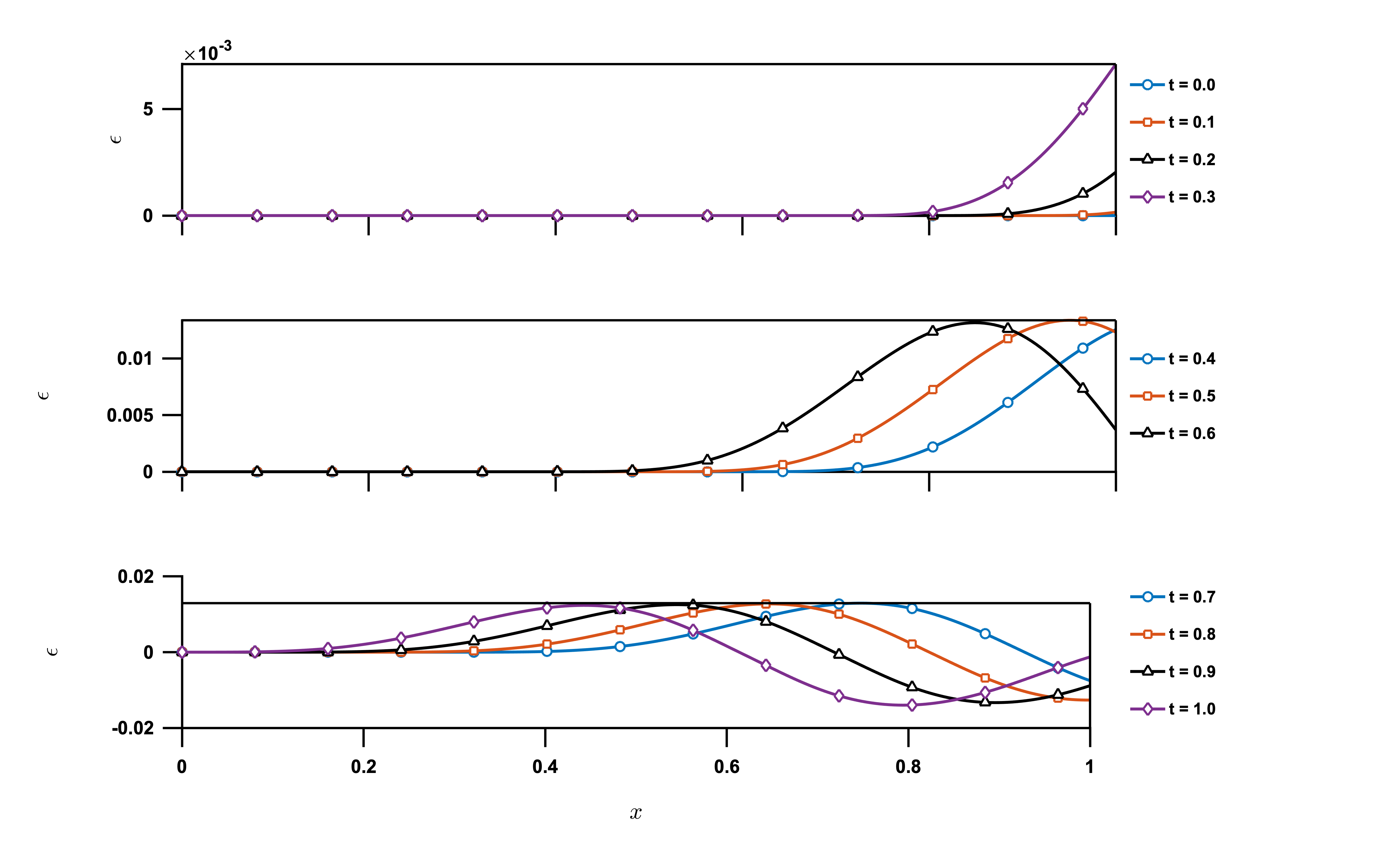}
        \caption{$\epsilon$}
        \label{fig:strain_beta5}
    \end{subfigure}  
    \caption{The stress $\sigma$ and strain $\epsilon$ for the moderately nonlinear case with $\beta=5.0$ and $\alpha=1.5$.}
    \label{fig:stress_strain_beta5}
\end{figure}

\subsection{Nonlinear case with $b=10.0$}

\begin{figure}[H]
    \centering
    \begin{subfigure}[b]{0.48\textwidth}
        \centering
        \includegraphics[width=\textwidth, height=9cm]{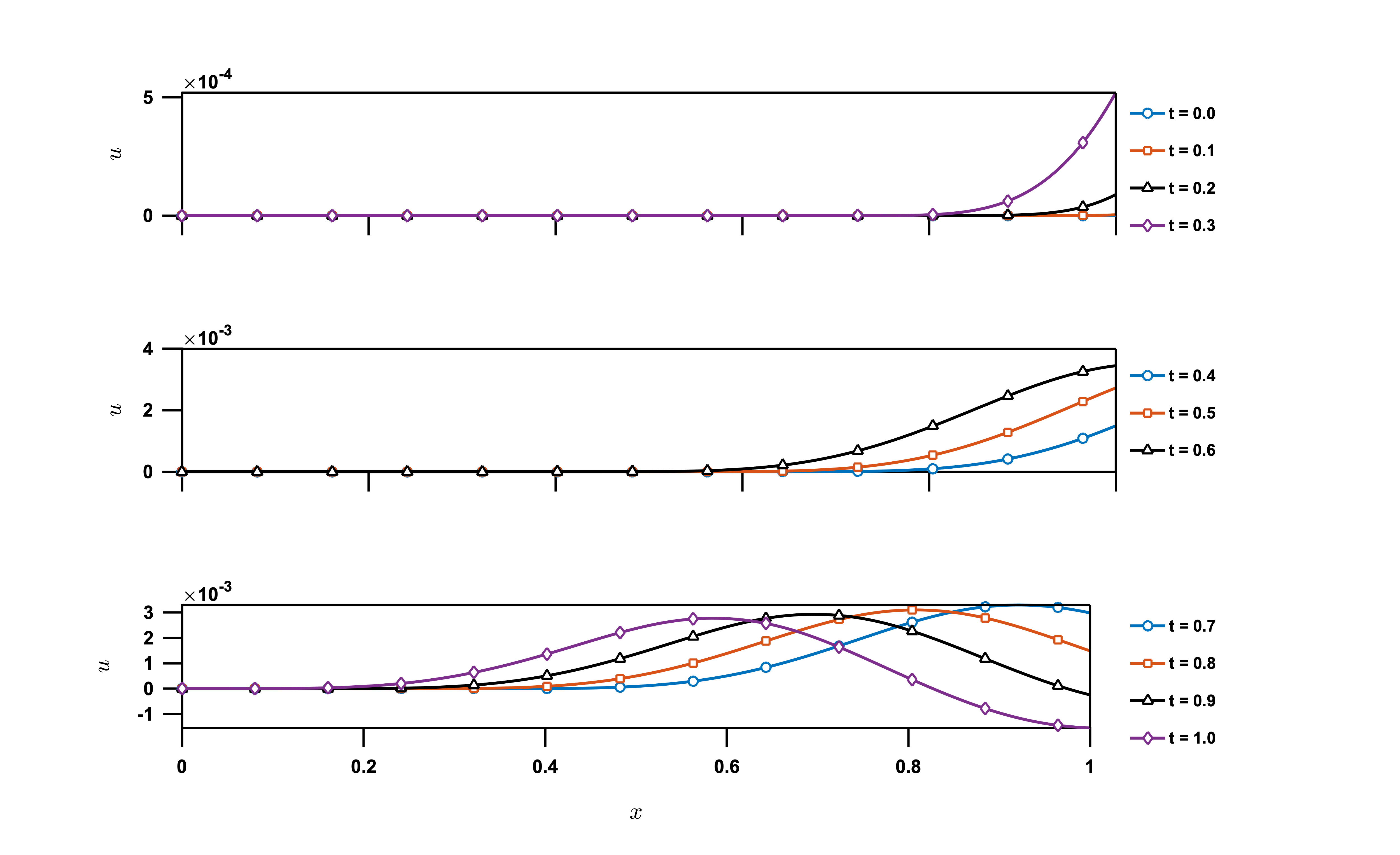}
        \caption{$u$}
        \label{fig:disp_b10}
    \end{subfigure}
    \hfill 
    \begin{subfigure}[b]{0.48\textwidth}
        \centering
        \includegraphics[width=\textwidth, height=9cm]{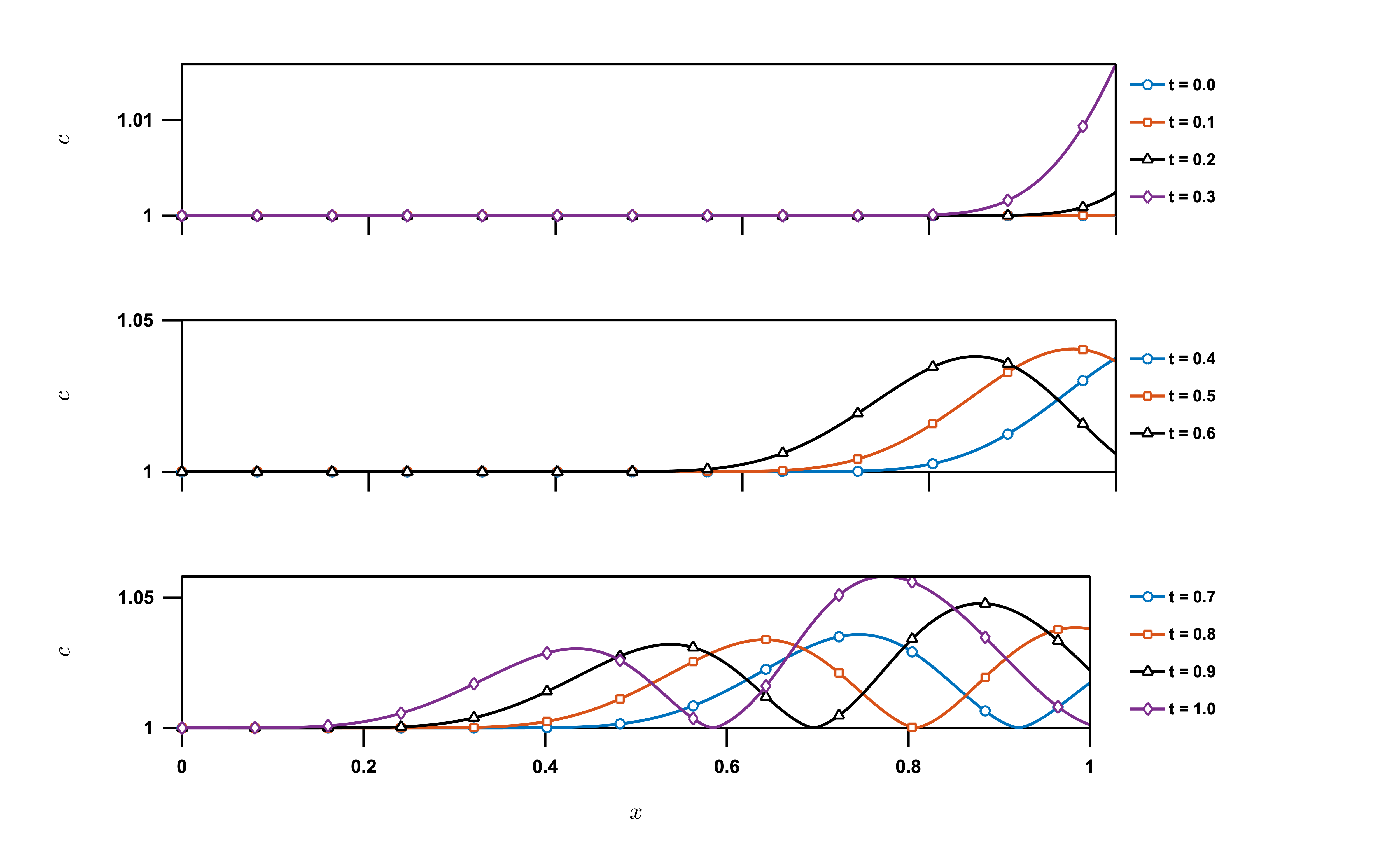}
        \caption{$c$}
        \label{fig:speed_b10}
    \end{subfigure}
    \caption{The displacement $u$ and wave speed $c$ for the nonlinear case with $b=10.0$ and $a=1.5$. Note the substantial increase in wave speed variation and the resulting wavefront steepening.}
    \label{fig:disp_speed_b10}
\end{figure}

Finally, we examine the highly nonlinear regime defined by $b = 10.0$ and $a = 1.5$. In this scenario, the constitutive nonlinearity exerts a dominant influence on the wave propagation physics. The wave speed $c$, as shown in Fig. \ref{fig:speed_b10}, demonstrates pronounced excursions from unity, reaching peak values of approximately $c \approx 1.04$. This significant stress-dependence results in extreme dispersive behavior. The stress $\sigma$ and strain $\epsilon$ profiles exhibit severe wavefront steepening as the pulse propagates leftward. Because the high-amplitude portions of the wave travel significantly faster than the low-amplitude tails, the trailing edge of the pulse sharpens dramatically, nearly leading to the formation of a shock discontinuity by $t=1.0$. This case clearly illustrates the transition from smooth wave translation to shock-dominated dynamics driven by the magnitude of the parameter $b$.

\begin{figure}[H]
    \centering
    \begin{subfigure}[b]{0.48\textwidth}
        \centering
        \includegraphics[width=\textwidth, height=9cm]{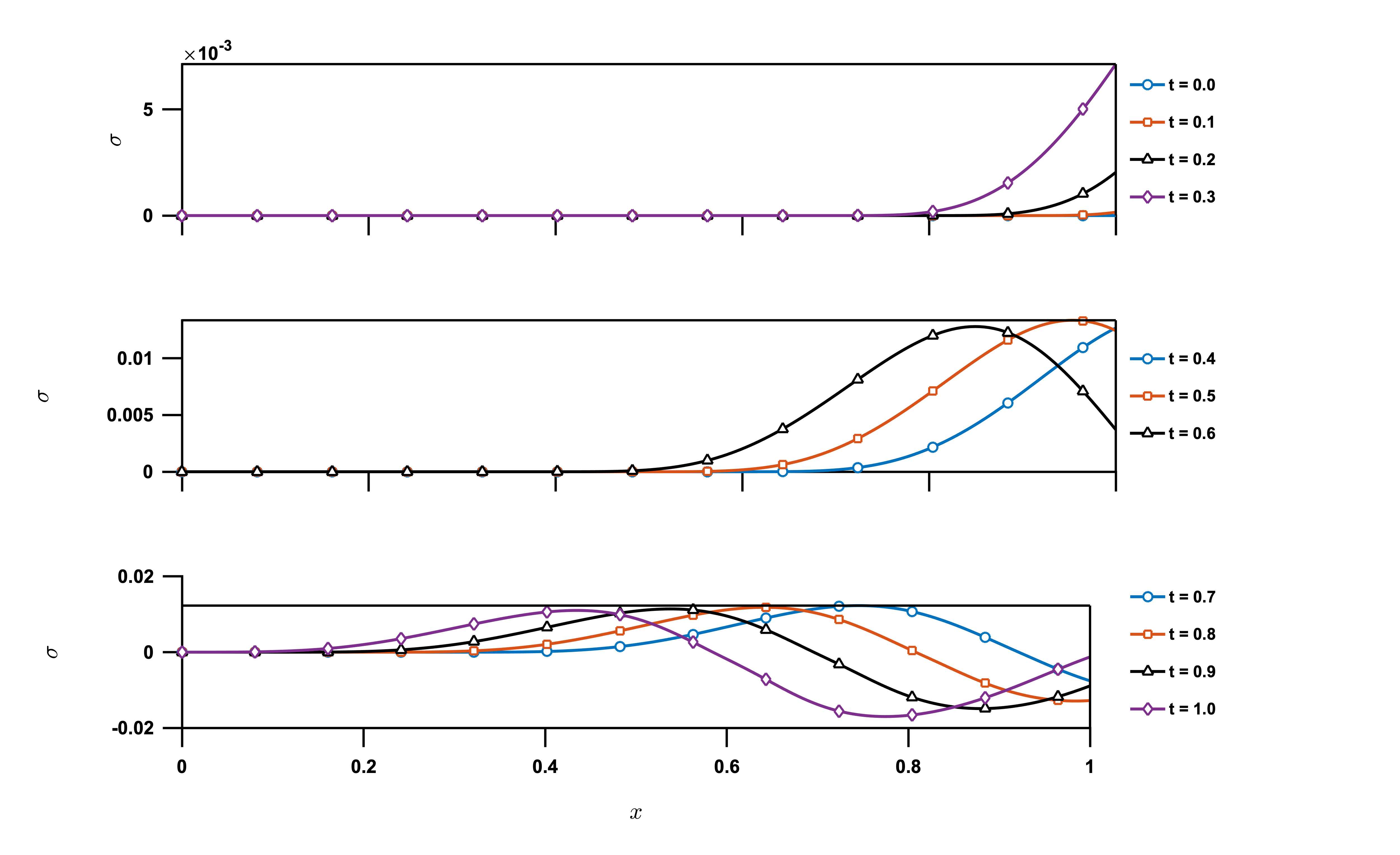}
        \caption{$\sigma$}
        \label{fig:sigma_b10}
    \end{subfigure}
    \hfill 
    \begin{subfigure}[b]{0.48\textwidth}
        \centering
        \includegraphics[width=\textwidth, height=9cm]{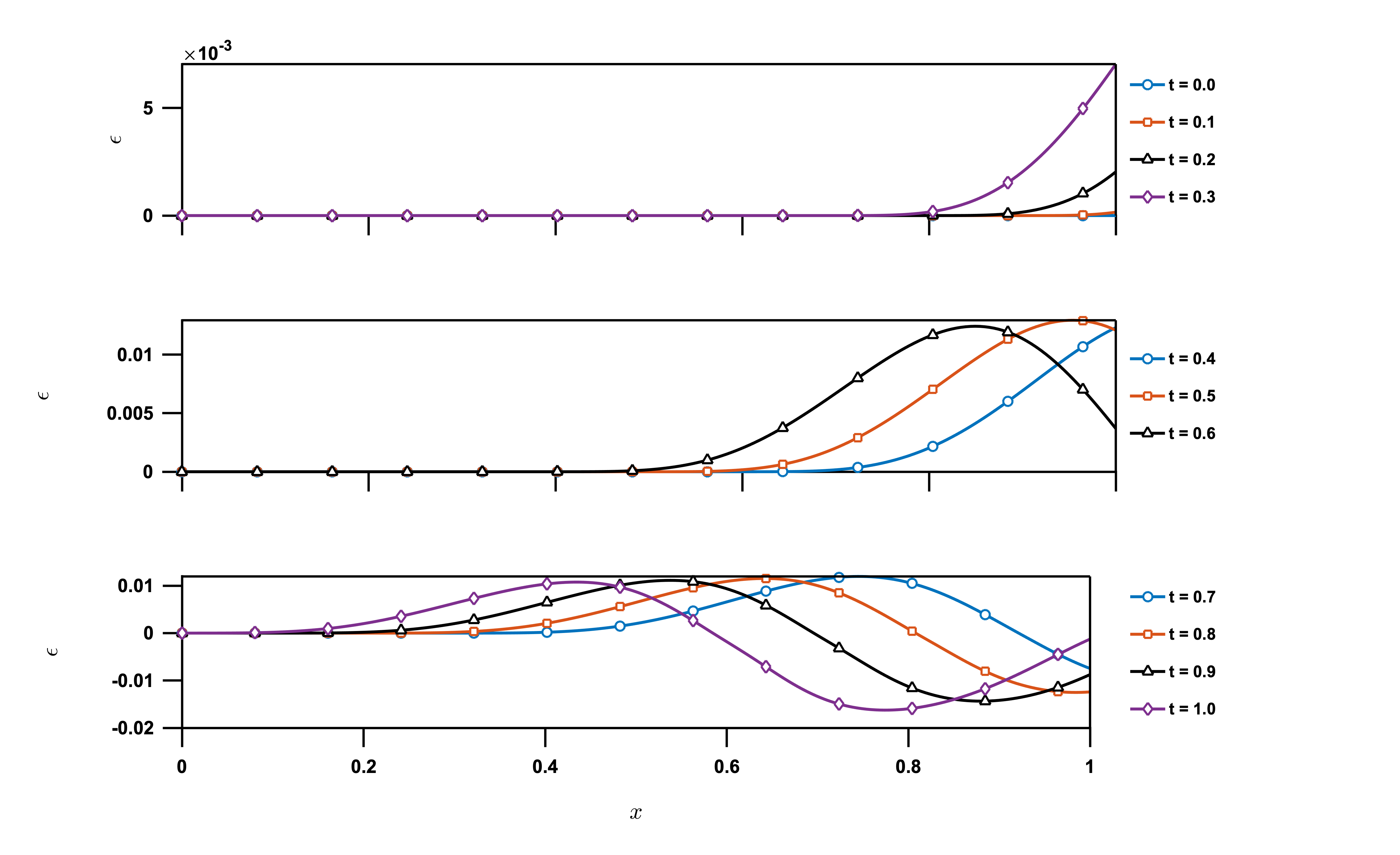}
        \caption{$\epsilon$}
        \label{fig:strain_b10}
    \end{subfigure}   
    \caption{The stress $\sigma$ and strain $\epsilon$ for the highly nonlinear case with $b=10.0$ and $a=1.5$.}
    \label{fig:stress_strain_b10}
\end{figure}

\subsection{Effect of parameter $a$ with fixed $b=1.0$}

This subsection investigates the sensitivity of the wave propagation to the exponent parameter $a$, while maintaining a constant coefficient $b=1.0$. We examine three distinct cases: $a=3.0$, $a=5.0$, and $a=10.0$. 

The results reveal a clear inverse relationship between the magnitude of $a$ and the observability of nonlinear phenomena in the small-strain regime ($\epsilon < 0.02$). Mathematically, since the nonlinear contribution scales with the strain raised to the power of $a$ (i.e., $\sim \epsilon^a$), and given that $\epsilon \ll 1$, increasing $a$ significantly reduces the magnitude of the nonlinear term. 

\begin{figure}[H]
    \centering
    \begin{subfigure}[b]{0.48\textwidth}
        \centering
        \includegraphics[width=\textwidth, height=9cm]{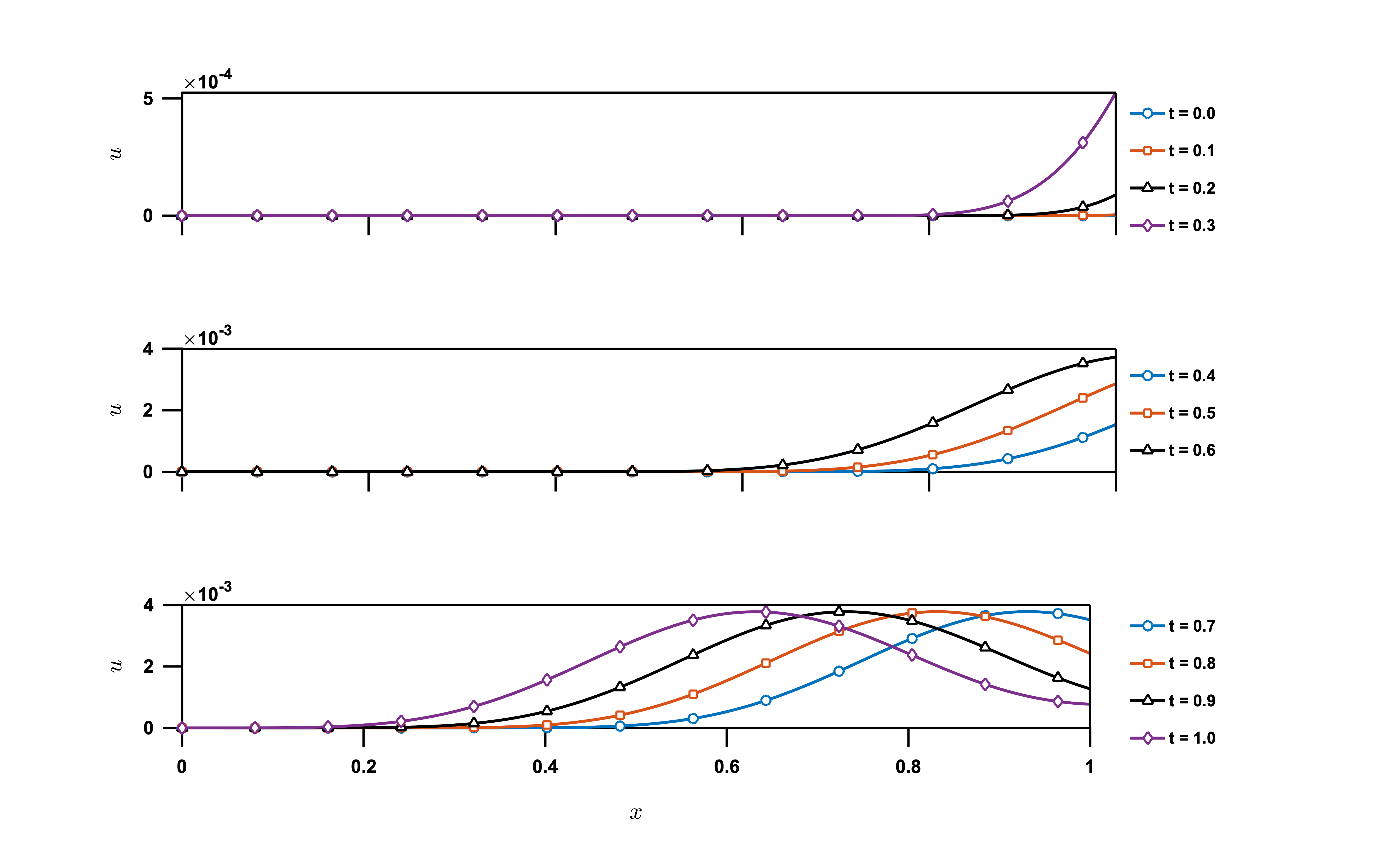}
        \caption{$u$}
        \label{fig:disp_a3}
    \end{subfigure}
    \hfill 
    \begin{subfigure}[b]{0.48\textwidth}
        \centering
        \includegraphics[width=\textwidth, height=9cm]{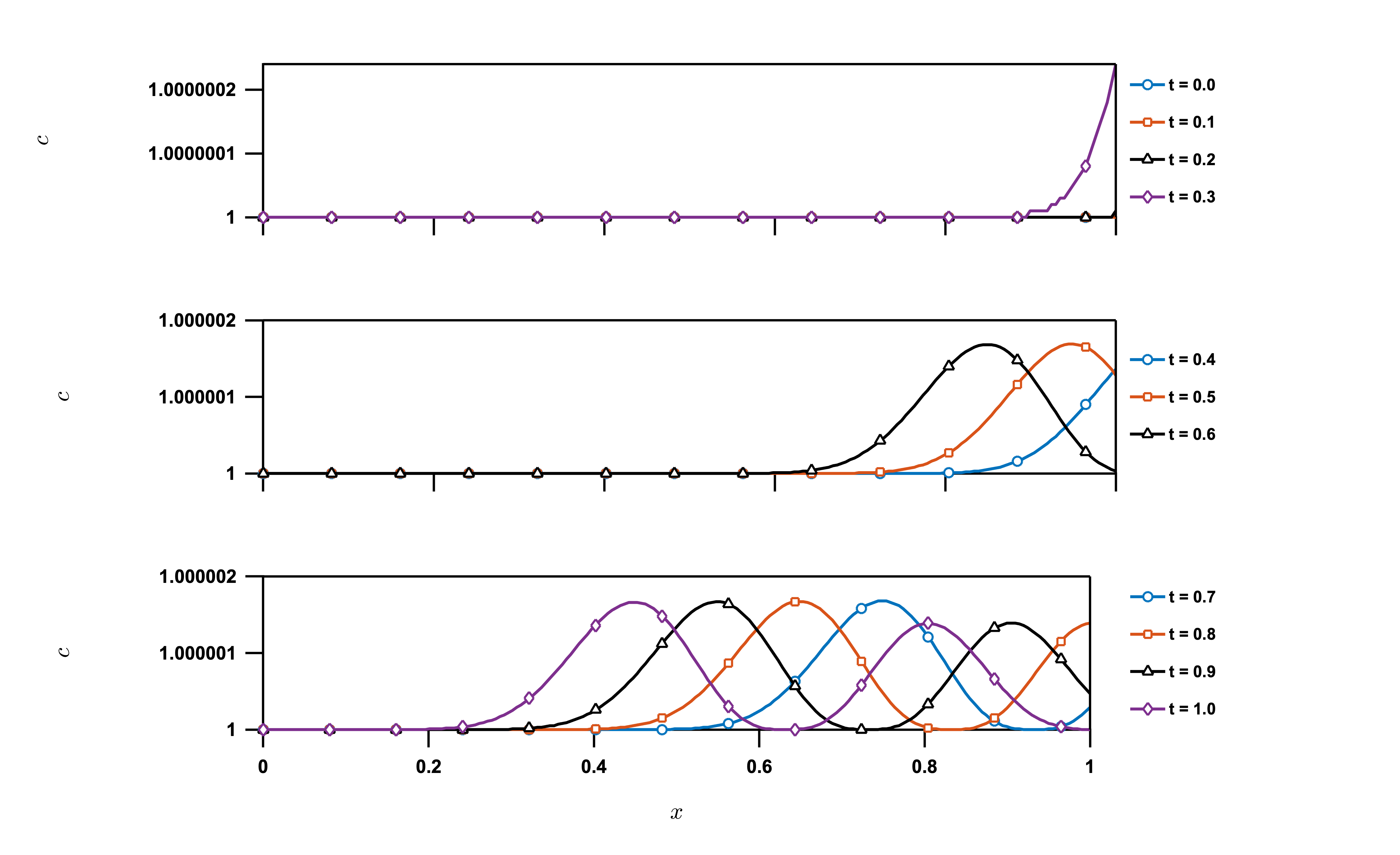}
        \caption{$c$}
        \label{fig:speed_a3}
    \end{subfigure}
    \caption{The displacement $u$ and wave speed $c$ for $b=1.0$ and $a=3.0$. Note the negligible deviation in wave speed.}
    \label{fig:disp_speed_a3}
\end{figure}

For the case of $a=3.0$, the wave speed $c$ (Fig. \ref{fig:speed_a3}) exhibits only minute deviations from unity (on the order of $10^{-6}$), rendering the system practically linear. As $a$ is increased further to $5.0$ and $10.0$, this linearization effect becomes even more pronounced. The nonlinear stiffness terms vanish rapidly, resulting in wave profiles for stress, strain, and displacement that maintain their shape without the steepening observed in lower-$a$ cases. Thus, the parameter $a$ effectively acts as a threshold controller; higher values confine nonlinear effects to regimes of significantly larger deformation, while lower values allow nonlinearity to manifest at small strain amplitudes.

\begin{figure}[H]
    \centering
    \begin{subfigure}[b]{0.48\textwidth}
        \centering
        \includegraphics[width=\textwidth, height=9cm]{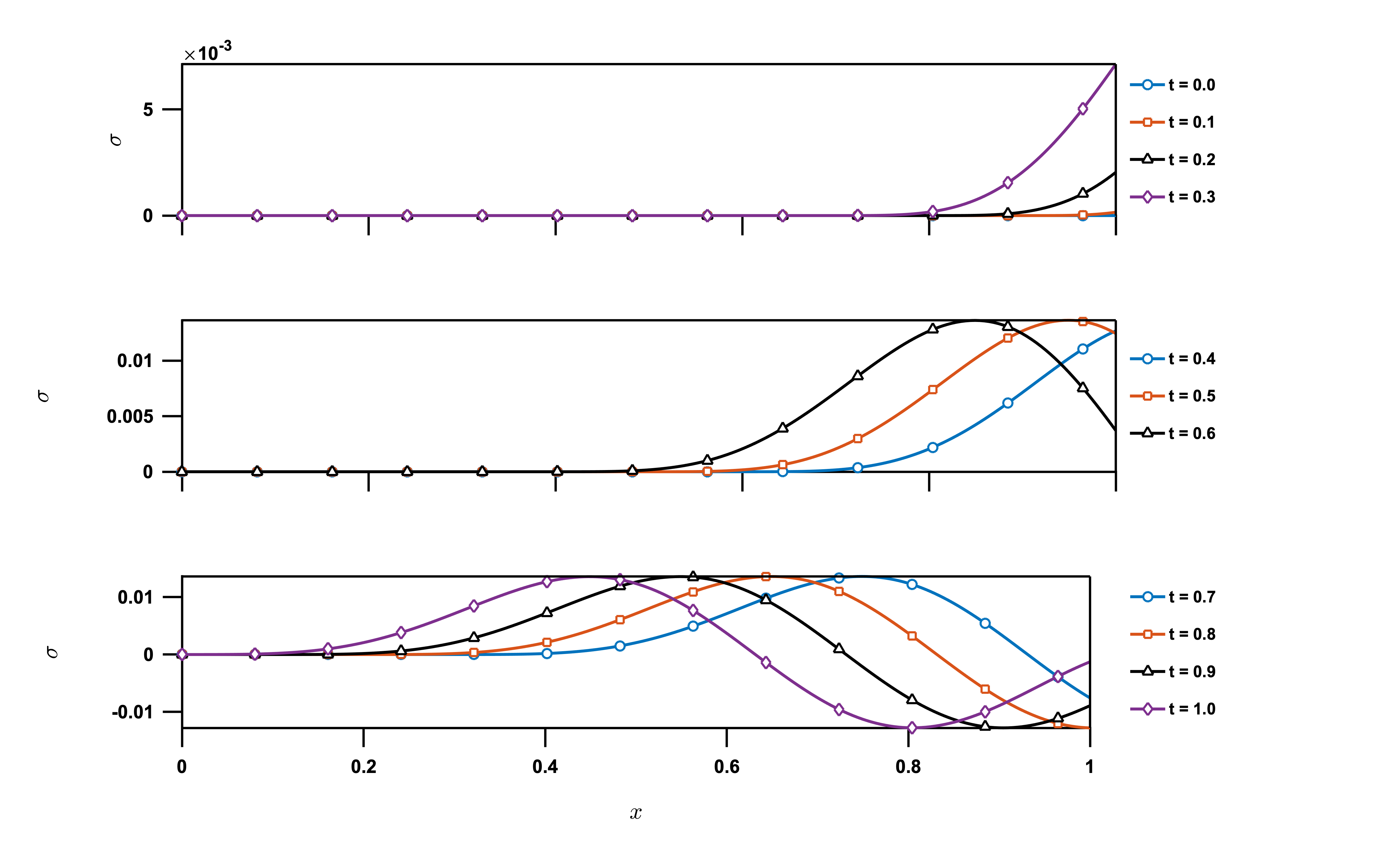}
        \caption{$\sigma$}
        \label{fig:sigma_a3}
    \end{subfigure}
    \hfill 
    \begin{subfigure}[b]{0.48\textwidth}
        \centering
        \includegraphics[width=\textwidth, height=9cm]{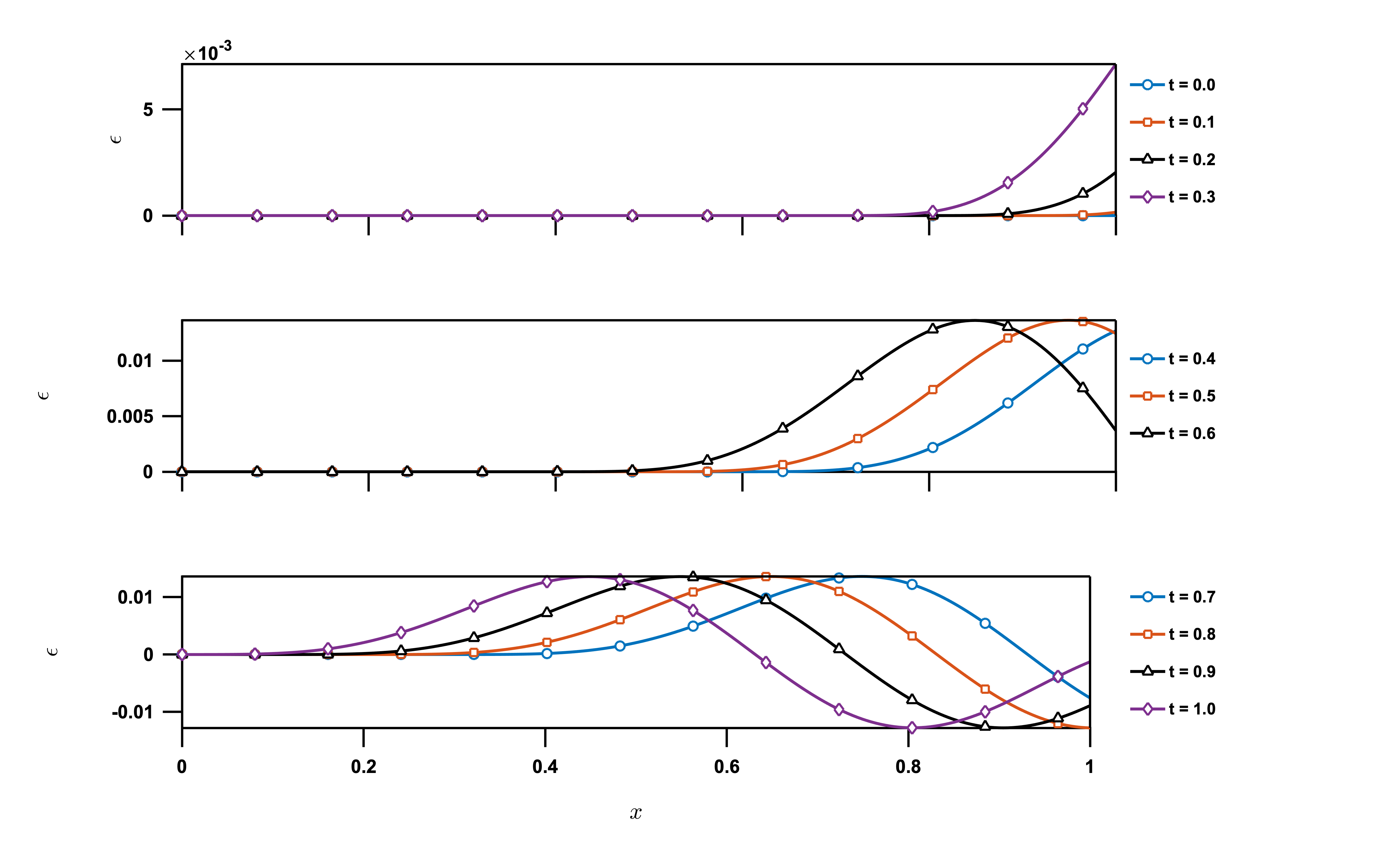}
        \caption{$\epsilon$}
        \label{fig:strain_a3}
    \end{subfigure}
    \caption{The stress $\sigma$ and strain $\epsilon$ for $b=1.0$ and $a=3.0$.}
    \label{fig:stress_strain_a3}
\end{figure}

\begin{figure}[H]
    \centering
    \begin{subfigure}[b]{0.48\textwidth}
        \centering
        \includegraphics[width=\textwidth, height=9cm]{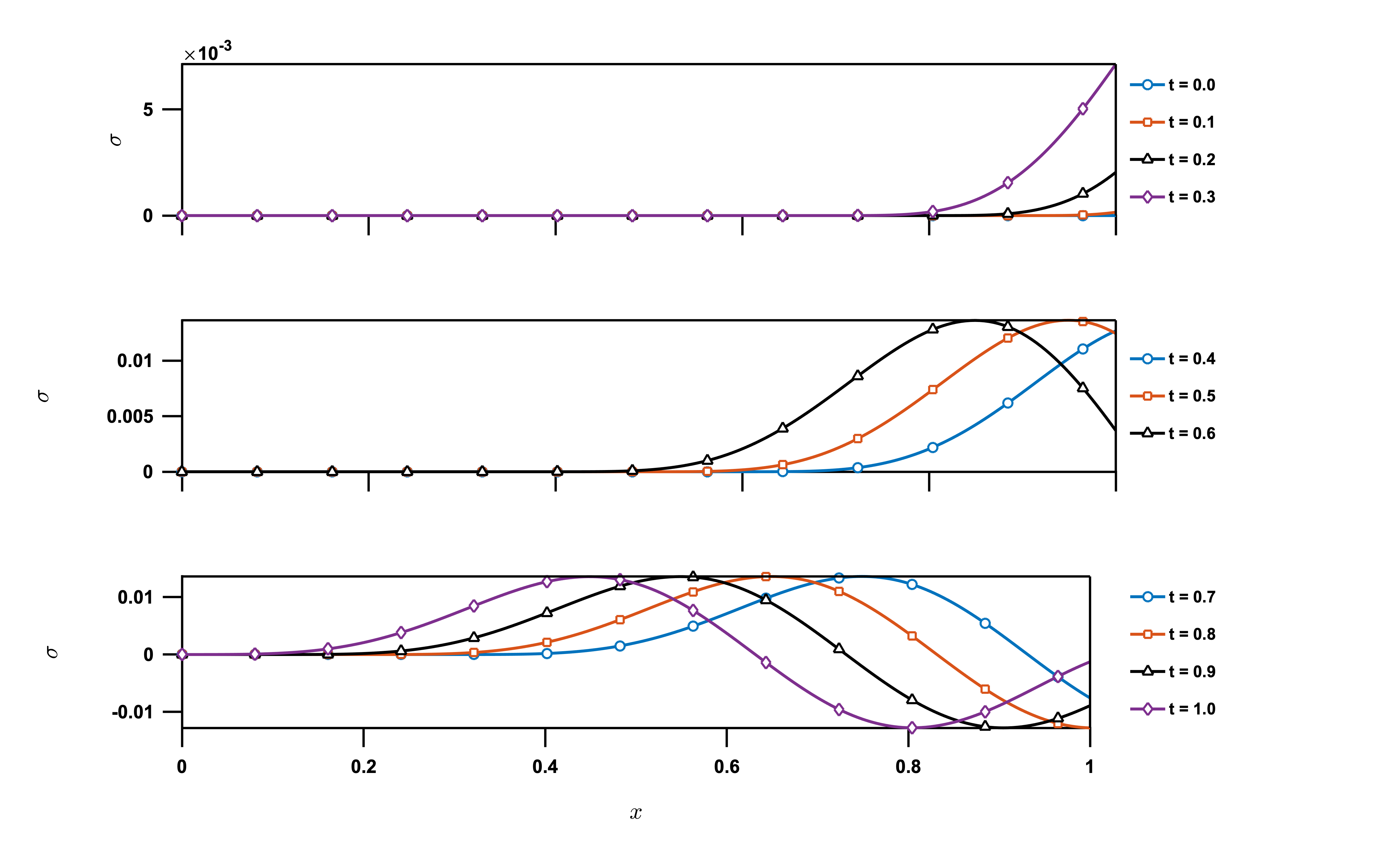}
        \caption{$\sigma$}
        \label{fig:sigma_a5}
    \end{subfigure}
    \hfill 
    \begin{subfigure}[b]{0.48\textwidth}
        \centering
        \includegraphics[width=\textwidth, height=9cm]{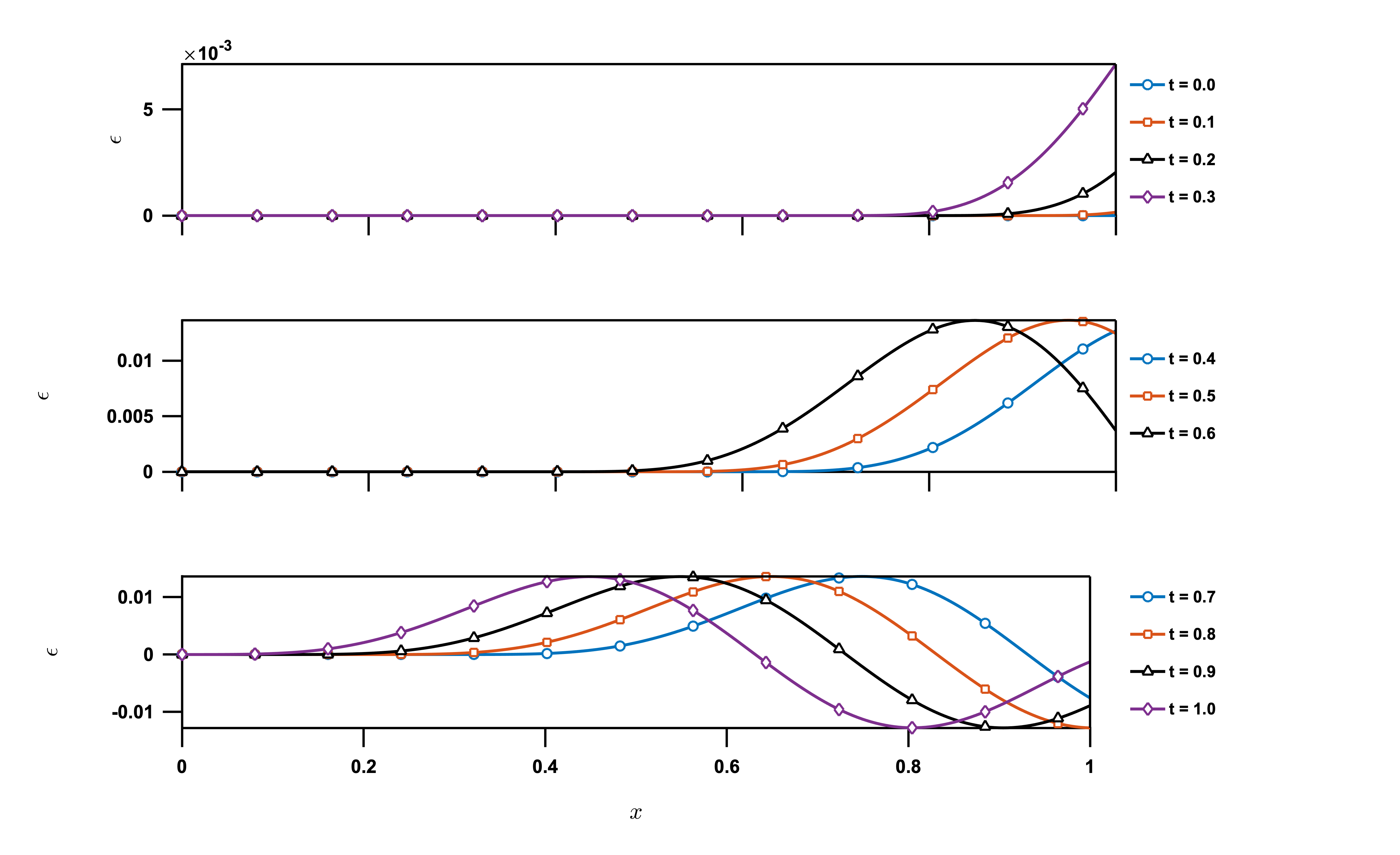}
        \caption{$\epsilon$}
        \label{fig:strain_a5}
    \end{subfigure}
    \caption{The stress $\sigma$ and strain $\epsilon$ for $b=1.0$ and $a=5.0$.}
    \label{fig:stress_strain_a5}
\end{figure}
\begin{figure}[H]
    \centering
    \begin{subfigure}[b]{0.48\textwidth}
        \centering
        \includegraphics[width=\textwidth, height=9cm]{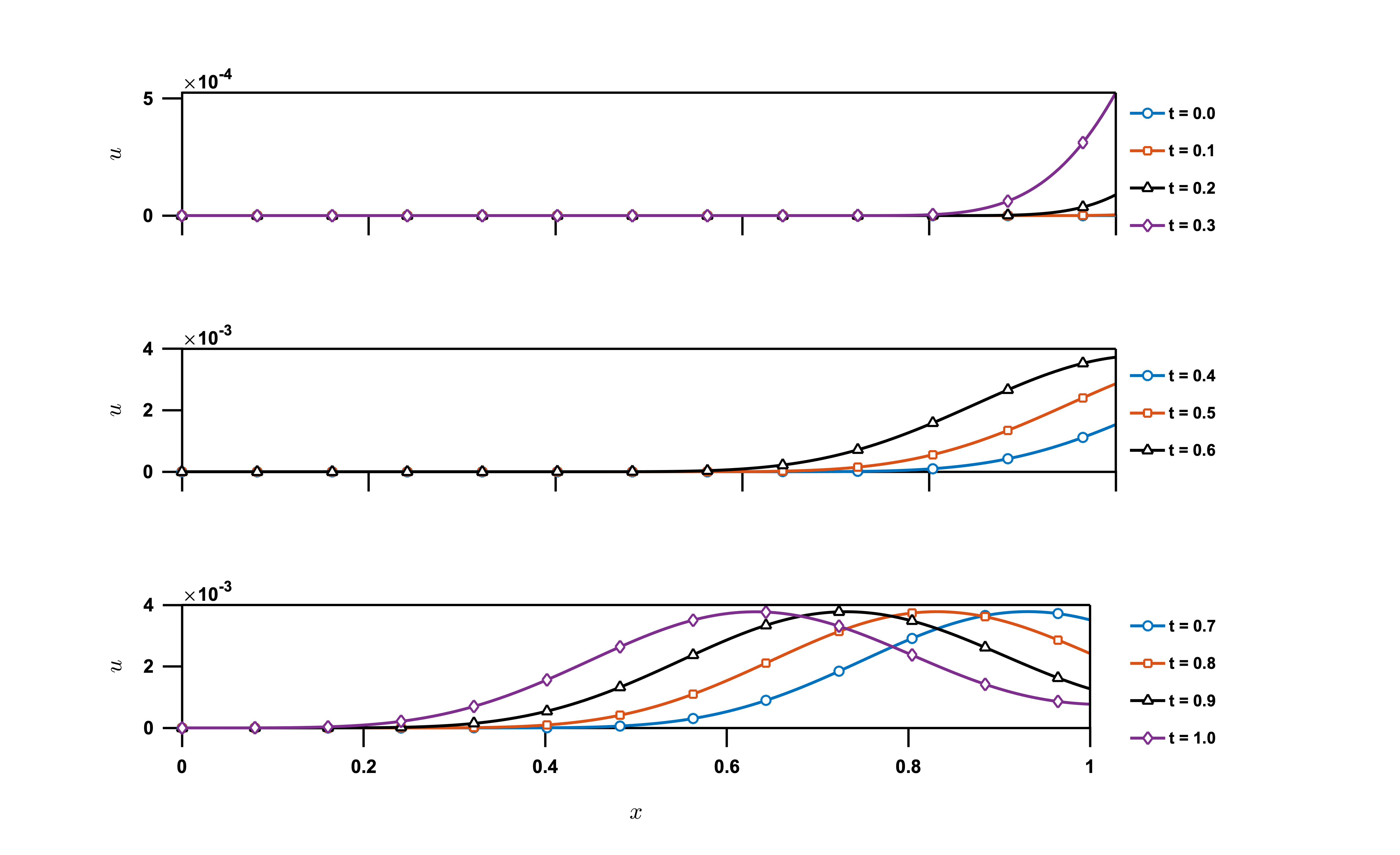}
        \caption{$u$}
        \label{fig:disp_a5}
    \end{subfigure}
    \hfill 
    \begin{subfigure}[b]{0.48\textwidth}
        \centering
        \includegraphics[width=\textwidth, height=9cm]{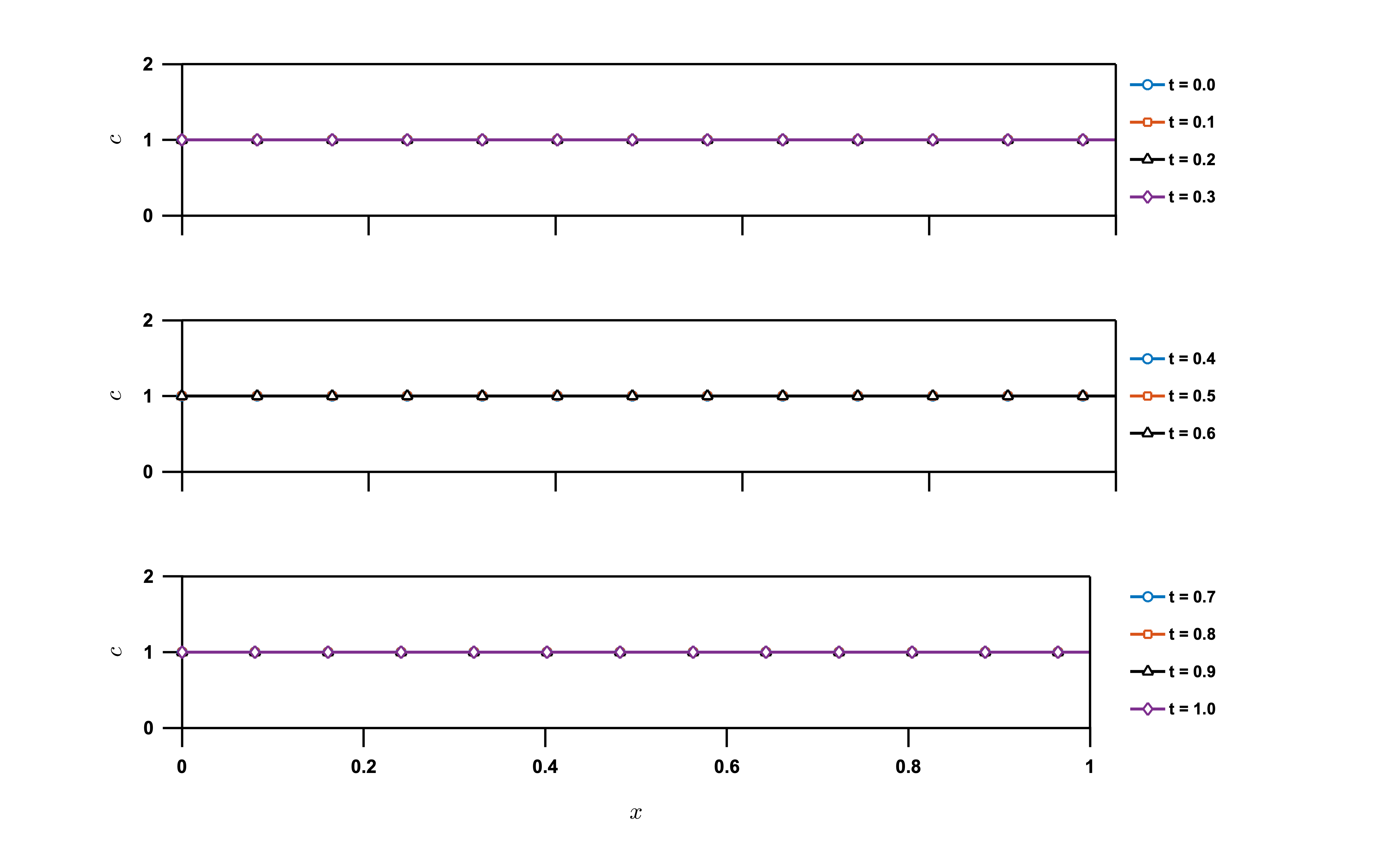}
        \caption{$c$}
        \label{fig:speed_a5}
    \end{subfigure}
    \caption{The displacement $u$ and wave speed $c$ for $b=1.0$ and $a=5.0$.}
    \label{fig:disp_speed_a5}
\end{figure}

\begin{figure}[H]
    \centering
    \begin{subfigure}[b]{0.48\textwidth}
        \centering
        \includegraphics[width=\textwidth, height=9cm]{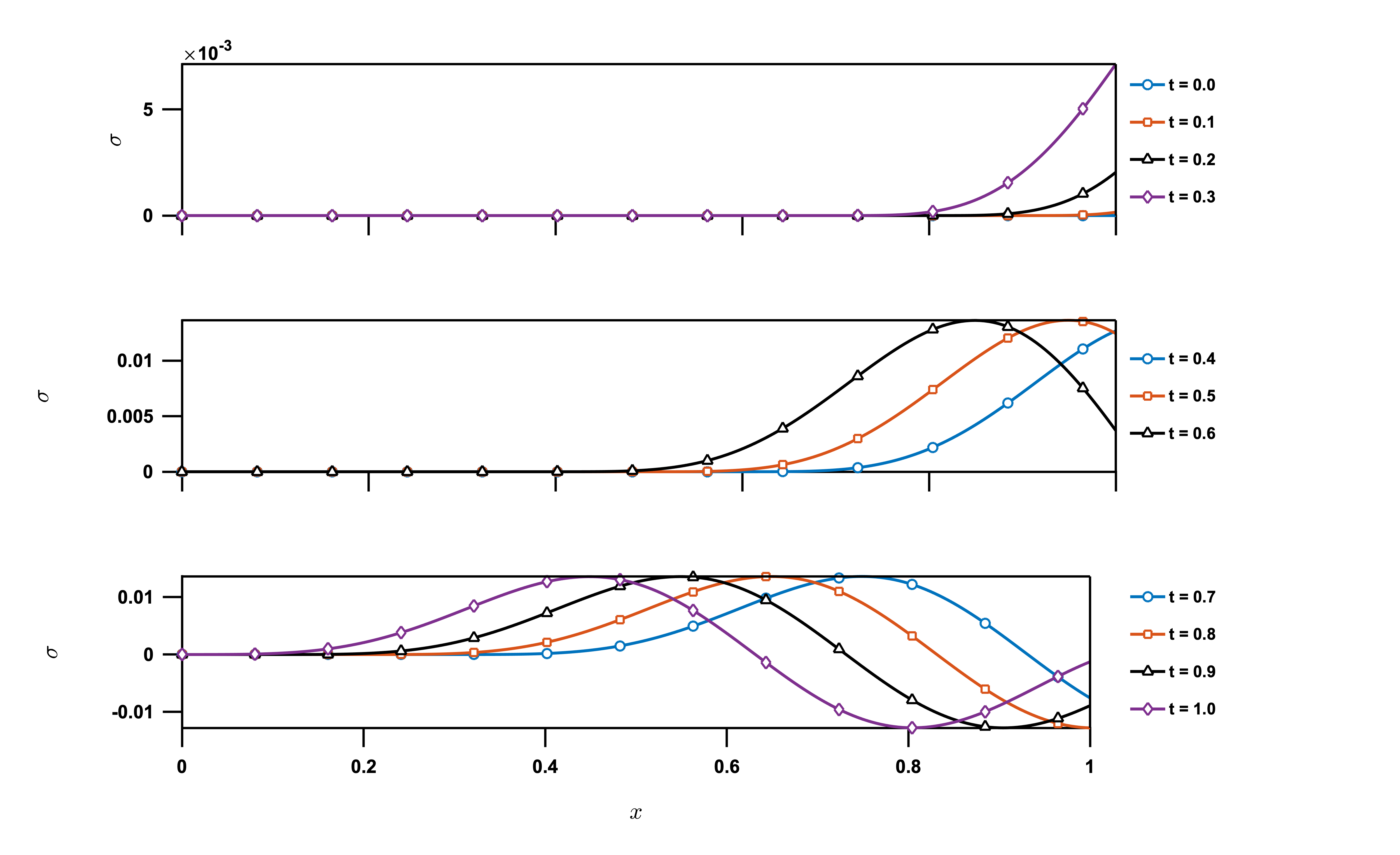}
        \caption{$\sigma$}
        \label{fig:sigma_a10}
    \end{subfigure}
    \hfill 
    \begin{subfigure}[b]{0.48\textwidth}
        \centering
        \includegraphics[width=\textwidth, height=9cm]{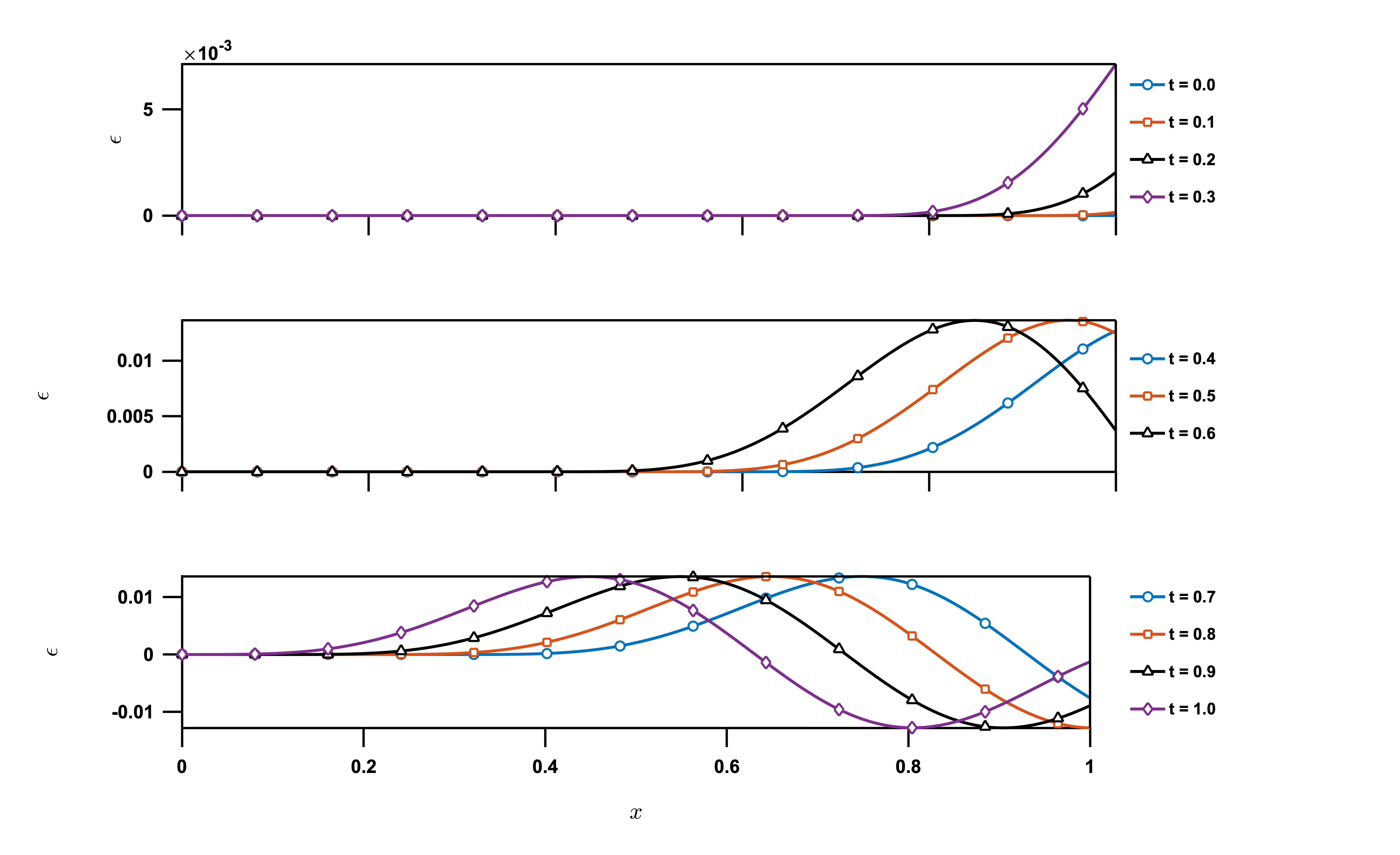}
        \caption{$\epsilon$}
        \label{fig:strain_a10}
    \end{subfigure}
    \caption{The stress $\sigma$ and strain $\epsilon$ for $b=1.0$ and $a=10.0$.}
    \label{fig:stress_strain_a10}
\end{figure}
\begin{figure}[H]
    \centering
    \begin{subfigure}[b]{0.48\textwidth}
        \centering
        \includegraphics[width=\textwidth, height=9cm]{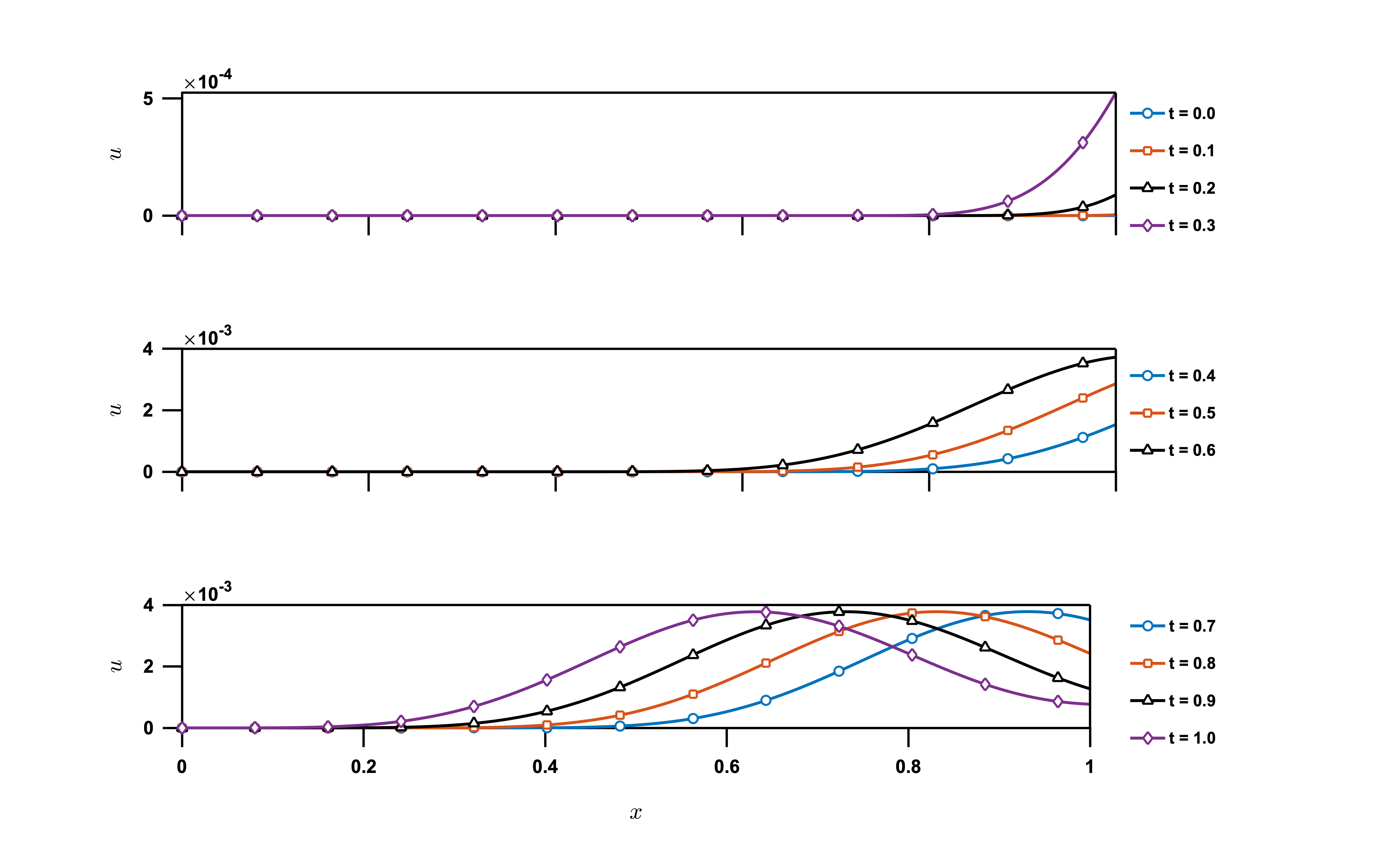}
        \caption{$u$}
        \label{fig:disp_a10}
    \end{subfigure}
    \hfill 
    \begin{subfigure}[b]{0.48\textwidth}
        \centering
        \includegraphics[width=\textwidth, height=9cm]{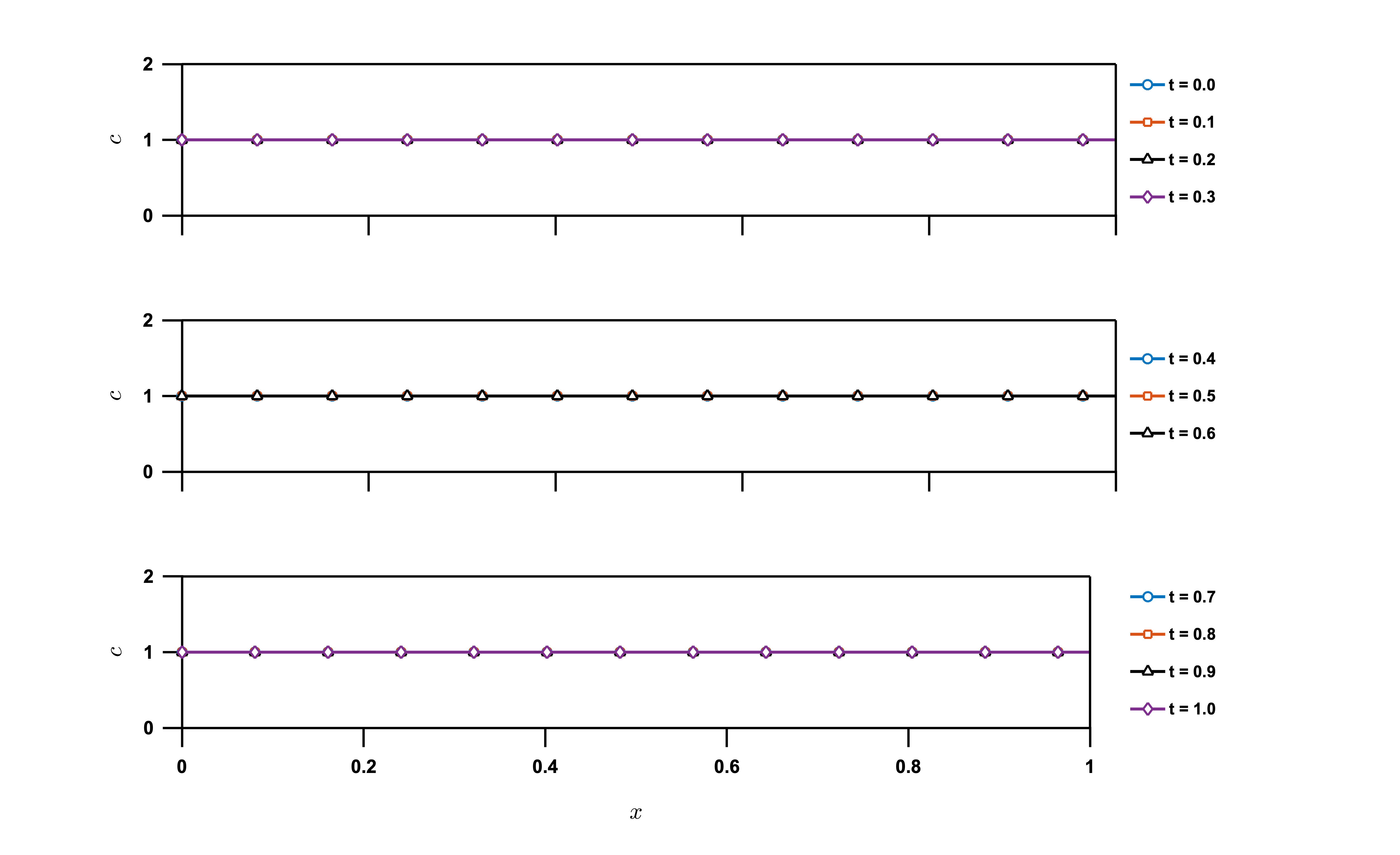}
        \caption{$c$}
        \label{fig:speed_a10}
    \end{subfigure}
    \caption{The displacement $u$ and wave speed $c$ for $b=1.0$ and $a=10.0$. The response is effectively linear due to the high exponent.}
    \label{fig:disp_speed_a10}
\end{figure}

\section{Conclusion} \label{sec:conclusion}
In this work, we presented a comprehensive computational study of wave propagation in geometrically linear elastic materials governed by algebraically nonlinear constitutive relations. By systematically varying the material parameters $b$ (magnitude of nonlinearity) and $a$ (exponent of nonlinearity), we established a clear mapping between the constitutive description and the resulting wave dynamics. The numerical framework was first validated against the linear elastic limit ($b=0$), recovering the expected shape-preserving wave translation and a constant unity wave speed.

The parametric study revealed distinct roles for the two constitutive parameters:
\begin{itemize}
    \item \textbf{Influence of parameter $b$:} The coefficient $b$ acts as the primary driver for the intensity of nonlinear effects. As $b$ was increased from 1.0 to 10.0, the wave speed exhibited increasingly significant spatial and temporal fluctuations. This stress-dependent wave velocity introduced dispersive mechanisms that led to the steepening of the wavefront. In the highly nonlinear regime ($b=10.0$), the differential wave speed between high- and low-amplitude regions resulted in the formation of sharp gradients, indicating a transition from smooth wave propagation to shock-dominated dynamics.
    
    \item \textbf{Influence of parameter $a$:} The exponent parameter $a$ controls the activation threshold of the nonlinearity. For the small-strain regimes investigated ($\epsilon < 0.02$), increasing the exponent $a$ effectively suppressed the nonlinear response. While lower values (e.g., $a=1.5$) allowed nonlinear features to manifest clearly, higher values (e.g., $a \ge 5.0$) rendered the nonlinear terms negligible, recovering an effectively linear response. This confirms that for power-law type materials, the exponent determines the strain threshold required to observe nonlinear wave physics.
\end{itemize}

In summary, the results demonstrate that the proposed numerical framework is robust and capable of capturing the complex interplay between constitutive nonlinearity and wave dispersion. The capability to accurately simulate the evolution from smooth pulses to shock discontinuities provides a valuable tool for characterizing advanced materials where the standard assumption of constant stiffness is insufficient. Future work will extend this formulation to multi-dimensional domains and explore the impact of implicit time-integration schemes on energy conservation in the presence of strong shocks.

\section*{Acknoledgement}
The work of SMM is supported by the National Science Foundation under Grant No.\ 2316905.

\bibliographystyle{plain}
\bibliography{references}

\end{document}